\documentclass[12pt,leqno]{compositio}

\usepackage{amsmath,amssymb,amscd,latexsym}

\usepackage[all]{xy}

\newcommand{\Br}{{\mathrm{Br}}}

\newcommand{\diag}{{\mathrm{diag}}}
\newcommand{\Frob}{{\mathrm{Frob}}}
\newcommand{\Gal}{{\mathrm{Gal}}}

\newcommand{\A}{\mathbb{A}}
       
\newcommand{\im}{{\mathrm{im}}}

\newcommand{\PGL}{{\mathrm{PGL}}}
\newcommand{\Sl}{{\mathrm{Sl}}}
\newcommand{\SO}{{\mathrm{SO}}}
\renewcommand{\O}{{\mathrm{O}}}
\newcommand{\Sp}{{\mathrm{Sp}}}
\newcommand{\Smo}{{\mathbb{S}}}
\newcommand{\SU}{{\mathrm{SU}}}
\newcommand{\T}{{\mathrm{T}}}
\newcommand{\U}{{\mathrm{U}}}

\newcommand{\Gl}{{\mathrm{Gl}}}

\newcommand{\Res}{{\mathrm{Res}}}

\newcommand{\Gm}{\mathbb{G}_m}

\newcommand{\F}{\mathbb{F}}

\newcommand{\G}{\mathbb{G}}

\newcommand{\C}{\mathbb{C}}

\newcommand{\HH}{{\mathrm{H}}}

\newcommand{\Nr}{{\mathrm{N}}}

\newcommand{\Q}{{\mathbb{Q}}}

\newcommand{\R}{{\mathbb{R}}}

\newcommand{\Z}{{\mathbb{Z}}}

\renewcommand{\dim}{\mathrm{dim}}

\newcommand{\inv}{\mathrm{inv}}

\newcommand{\id}{\mathrm{id}}

\renewcommand{\ker}{\mathrm{ker}}

\newcommand{\rk}{\mathrm{rk}}

\newcommand{\Spec}{\mathrm{Spec}}

\newcommand{\Aut}{\mathrm{Aut}\,}

\newcommand{\End}{\mathrm{End}\,}

\newcommand{\M}{\mathrm{M}}

\newcommand{\Rf}{\mathcal{R}}

\newcommand{\Srm}{\mathrm{S}}

\newcommand{\Df}{\mathcal{D}}

\newcommand{\tr}{\mathrm{tr}}

\newcommand{\Oh}{{\mathcal O}}

\newcommand{\OSfin}{{\mathcal O}_{k,{\mathrm S}^{\mathrm{ fin}}}}

\newcommand{\p}{{\mathfrak{p}}}

\newcommand{\qf}{{\mathfrak{q}}}

\renewcommand{\P}{{\mathfrak{P}}}

\newcommand{\isoto}{\stackrel{\sim}{\longrightarrow}}

\newcommand{\st}{\, | \,}

\newtheorem{theorem}{Theorem}

\newtheorem{prop}[theorem]{Proposition}

\newtheorem{cor}[theorem]{Corollary}

\newtheorem{example}[theorem]{Example}

\newtheorem{remark}[theorem]{Remark}

\newenvironment{proofof}{\noindent {\bf Proof of}}{\mbox{}\hfill$\Box$}

\begin{document}

\title[Arithmetically defined dense subgroups of\\ Morava stabilizer groups]{Arithmetically defined dense subgroups of \\ Morava stabilizer groups}

\author{Niko Naumann}

\email{niko.naumann@mathematik.uni-regensburg.de}

\address{NWF I- Mathematik\\ Universit\"at Regensburg\\93040 Regensburg}



\begin{abstract} For every prime $p$ and integer $n\ge 3$ we
  explicitly
construct an abelian variety $A/\F_{p^n}$ of dimension $n$ such that
for a suitable prime $l$ the group of quasi-isogenies of $A/\F_{p^n}$
of
$l$-power degree is canonically a dense subgroup of the $n$-th Morava 
stabilizer group at $p$.
We also give a variant of this result taking into account a
polarization. 
This is motivated by the recent construction of topological automorphic forms which generalizes topological modular forms \cite{topautforms}.\\
For this, we prove some arithmetic results of independent interest: a result about approximation of local units in maximal orders of global skew-fields which also gives a precise
 solution to the problem of extending automorphisms of the 
$p$-divisible group of a simple abelian variety over a finite field
to quasi-isogenies of the abelian variety of degree divisible by as
few primes as possible. 

\end{abstract}

\maketitle

\section{Introduction}\label{intro}
One of the most fruitful ways of studying the stable homotopy category
is the chromatic approach: after localizing, in the sense of
Bousfield, at a prime $p$, one is left with an infinite hierarchy
of primes corresponding to the Morava $K$-theories $K(n)$, $n\ge 0$,
see
\cite{ravenelbook}. The successive layers in the resulting filtration are
the $K(n)$-local categories \cite{hoveystrickland} the structure of
which is governed to a large extend by (the continuous cohomology of)
the $n$-th Morava stabilizer group $\Smo_n$, i.e. the automorphism group
of 
the one-dimensional commutative formal group of height $n$ over
$\overline{\F}_p$. A fundamental problem in this context is to
generalize the fibration

\[ L_{K(1)}S^0\longrightarrow E_1^{hF}\longrightarrow E_1^{hF}, \]

c.f. the introduction of \cite{hennetc}, to a resolution of the $K(n)$-local sphere for $n\ge 2$. Substantial progress on this problem for $n=2$ and in many other cases as well has been achieved by clever use of
homological
algebra for $\Smo_n$-modules \cite{hennetc},\cite{henn2}. Recently, pursuing a
question of M. Mahowald and C. Rezk, M. Behrens was able to give a modular interpretation
of one such resolution in the case $n=2$ \cite{behrens}.\\
A basic observation is that $\Smo_2$ is the automorphism group of the
$p$-divisible group of a super-singular elliptic curve $E$ over a
finite
field $k$. Hence it seemed plausible, and was established in {\em
  loc. cit.},
that the morphisms in a resolution of a spectrum closely related to $L_{K(2)} S^0$ should have a description in terms
of
suitable endomorphisms of $E$. A key result for seeing this was to
observe that for suitable primes $l$

\begin{equation}\label{denseintro}
\left( \End_k(E)\left[ \frac{1}{l} \right] \right)^* \subseteq \Smo_2
\end{equation}

is a dense subgroup \cite[Theorem 0.1]{bl}.\\
One of our main results, Theorem \ref{modular1}, is the direct
generalization of (\ref{denseintro}) to arbitrary chromatic level
$n\ge 3$ in which $E$ is replaced by an abelian variety of dimension
$n$ which is known to be the minimal dimension possible.\\
In Corollary \ref{unitaryapplication} we give a variant of the arithmetic result underlying Theorem \ref{modular1}
in which on the left-hand-side of (\ref{denseintro}) we only allow
endomorphisms which are unitary with respect to a given
Rosati-involution. The motivation for this stems from recent work of M. Behrens and T. Lawson \cite{topautforms} bringing the arithmetic of suitable (derived) Shimura varieties to bear on homotopy theoretic
problems of arbitrary chromatic level, generalizing the role of
topological modular forms for problems of chromatic level at most
two, c.f. \cite[Theorem 15.2.1]{topautforms}.\\
This paper is organized as follows.\\
In subsection \ref{torsor} we record a well-known result about generically trivial torsors for later reference.
In subsection \ref{geo} we determine the structure of certain naturally occurring integral models for forms of $\Sl_d$, see
Theorems \ref{smooth} and \ref{geoofUandSU}.\\
As a first application, in section \ref{approxlocal}, we consider the problem
of approximating local units of maximal orders in finite-dimensional 
skew-fields over $\Q$ (carrying a positive involution of the second kind)
by global (unitary) units with as few denominators as possible. This is
naturally an approximation problem for specific integral models of the general linear (certain unitary) group(s) and will be reduced to a similar problem 
for $\G_m$ (a specific integral model $\T'$ of a one-dimensional anisotropic
torus) in Theorems \ref{approx} and \ref{approxinv}.\\
In section \ref{commutative}
we can solve the approximation problem for $\G_m$ using class field theory and settle a special case for $\T'$, see Theorems \ref{commapprox} and \ref{approxtorus}.\\
In subsection \ref{extendaut} we explain the application of the
results   
obtained so far to the following problem: given a simple abelian
variety $A$ over a finite field one would like to extend an
automorphism of the $p$-divisible group $A[p^{\infty}]$ of $A$ to a 
quasi-isogeny of $A$ the degree of which should be divisible by as few
primes as possible. Subsection \ref{qisogdense} contains the proof of Theorem \ref{modular1} reviewed above.\\

\begin{acknowledgements} 

I would like to thank U. Jannsen and A. Schmidt for useful discussion concerning subsection \ref{cad-1} and M. Behrens and T. Lawson for pointing out to me the
abelian varieties used in subsection \ref{qisogdense}. Furthermore, I thank M. Behrens for pointing out a gap in the first version of this paper and J. Heinloth for his help to fix it. Finally, I am grateful to one
referee for suggesting stylistic improvements and to the other for a very careful report which led to substantial simplifications.

\end{acknowledgements}

\section{Arithmetic}\label{arithmetic}

\subsection{Triviality of torsors}\label{torsor}

The following result is well-known but we wish to state it in the form most suitable
for later references.

\begin{prop}\label{kertrivial} Let $k$ be a number-field with ring of integers $\Oh_k$, $\emptyset\neq U\subseteq\Spec(\Oh_k)$ an open sub-scheme and $G/U$ an affine smooth group-scheme. Assume that $G/U$ has connected fibers, that the generic fiber $G_k=G\times_U\Spec(k)$ is $k$-simple, semi-simple and simply-connected and that there is a place $v$ of $k$ outside $U$ such that $G$ is isotropic at $v$.\\
Then the canonical map of pointed sets \[ \HH^1(U,G)\to\HH^1(\Spec(k),G)\]
has trivial kernel.
\end{prop}
\begin{proof}
We use a result of A. Nisnevich, see \cite[Th\'eor\`eme 5.1]{gille}. Since $G/U$
has connected fibers we have $\HH^1_{fppf}(\hat{\Oh_x},G)=0$ in the notation
of {\em loc. cit.} It is thus sufficient to see that for every finite set $\Sigma$ of closed points of $U$ we have
\begin{equation}\label{einelement}
\mid \left( \prod\limits_{\p\in\Sigma}G(\Oh_{k,\p}) \setminus G(k_{\p})\right) /G(U-\Sigma)\mid=1.
\end{equation}

Here $\Oh_{k,\p}$ is the completion of $\Oh_k$ at $\p$ and $k_{\p}$ is the
field of fractions of $\Oh_{k,\p}$. The proof of (\ref{einelement}) using
strong approximation is very similar to the proof of Proposition \ref{density}
and is therefore left to the reader.
\end{proof}

\subsection{The geometry of some groups}\label{geo}

In this subsection we consider forms of
$\Sl_d$. These can be described in terms of skew-fields (with involution). The choice of a maximal order in the skew-field determines an integral model of the algebraic group and we will study the geometry of these group-schemes. The referee points out that some of these results, notably Theorem \ref{smooth} and 
Theorem \ref{geoofUandSU}, part {\sl ii)} are part of Bruhat-Tits theory, c.f. \cite[\S 5]{bt2} and \cite[4.6.32]{bt1}.

\subsubsection{Type $A_{d-1}$}\label{geoad-1}\hfill

Let $D$ be a finite dimensional skew-field over $\Q$ and $\Oh\subseteq D$ a maximal order \cite[Kapitel VI]{deu}. The center of $D$, denoted $k$, is a number field and we denote by $d$
the reduced dimension of $D$, i.e. $\dim_k\, D=d^2$. We denote by $\Oh_k\subseteq k$
the ring of integers and note that $k\cap\Oh=\Oh_k$ as an immediate consequence
of \cite[Kapitel VI,\S 11, Satz 7]{deu}.\\
Recall that $D$ is determined by its local invariants as follows \cite[Section 1.5.1]{pr}. Writing $\Sigma_k$ for the set of places of $k$, for every $v\in\Sigma_k$ there is a local invariant 
$\inv_v(D)\in\frac{1}{d}\Z/\Z\subseteq\Q/\Z$ and $\inv_v(D)=0$ for almost
all $v$. For a given place $v$, we denote by $k_v$ the completion of $k$
with respect to $v$. Then $D_v:=D\otimes_k k_v$ is a central simple $k_v$-algebra 
which determines a class $[D_v]\in\Br(k_v)$ in the Brauer group of $k_v$. There are
specific isomorphisms

\[ \tau_v:\Br(k_v)\isoto \left\{ \begin{array}{ccl} \Q/\Z & , & v\mbox{ finite}\\
\frac{1}{2}\Z/\Z & , & v\mbox{ real}\\ 
0 & , &   v \mbox{ complex} \end{array} \right. \]

such that $\inv_v(D)=\tau_v([D_v])$. Note that for every $v\in\Sigma_k$, $D_v$ is a skew-field if and
only if the order of $\inv_v(D)$ is $d$. \\
The group-valued functor $G$ on $\Oh_k$-algebras $R$

\[ G(R):=(\Oh\otimes_{\Oh_k} R)^*\]

is representable by an affine group scheme of finite type $G/\Spec(\Oh_k)$. The 
reduced norm induces a morphism of group-schemes $\Nr:G\longrightarrow\G_m$ over $\Spec(\Oh_k)$ and writing $G':=\ker(\Nr)$ gives an exact sequence of representable $fppf$-sheaves on $\Spec(\Oh_k)$

\begin{equation}\label{sequence}
1\longrightarrow G'\longrightarrow G\stackrel{\Nr}{\longrightarrow} \Gm \longrightarrow 1. 
\end{equation}

To see that $\Nr$ is $fppf$-surjective, note that the inclusion $R^*\subseteq(\Oh\otimes_{\Oh_k} R)^*$ defines a closed immersion $i: \G_m\longrightarrow G$ such that $\Nr\circ i$ is multiplication by
$d$ as can be checked on the generic fiber.

\begin{theorem}\label{smooth}
The groups $G$ and $G'$ are smooth with connected fibers over $\Spec(\Oh_k)$.
\end{theorem}

For the proof, we will need the following which might be compared with \cite[Expos\'e VI$_{B}$, Proposition 9.2 (xi)]{sga3} 

\begin{prop}\label{flat1}
Let $S$ be a scheme, $G,H/S$ group-schemes of finite presentation with affine fibers and $G/S$ flat and let $\phi:G\longrightarrow H$ be a morphism of $S$-group-schemes. Then the following are equivalent and imply that $H/S$ is flat:\\
i) $\phi$ is faithfully flat.\\
ii) $\phi$ is an epimorphism of $fppf$-sheaves.\\
iii) For every geometric point $\Spec(\Omega)\longrightarrow S$, $\phi_{\Omega}$ is an epimorphism of $fppf$-sheaves.
\end{prop}

\begin{proof} 
Since $\phi$ is of finite presentation, the implications {\sl i)}$\Rightarrow${\sl ii)}$ \Rightarrow ${\sl iii)} are trivial, so assume 
that {\sl iii)} holds true. Then, for every geometric point 
$\Spec(\Omega)\longrightarrow S$, the morphism of $\Omega$-Hopf algebras corresponding to $\phi_{\Omega}$ is injective, this follows from the existence of a $fppf$-local section of $\phi_{\Omega}$, and thus faithfully flat \cite[Theorem 4.1]{waterhousegroups}. This shows that $\phi$ is surjective and the fiber-wise
criterion for flatness \cite[Corollaire 11.3.11, $a)\Rightarrow b)$]{ega43} implies that $\phi$ and $H/S$ are flat.
\end{proof}

\begin{proofof} Theorem \ref{smooth}. To see that $G/\Spec(\Oh_k)$ is smooth we use the lifting criterion \cite[Remarques 17.1.2,i) et 17.5.4]{ega44}: if $A\longrightarrow A/I$ is the quotient of an Artinian
$\Oh_k$-algebra $A$ by an ideal $I\subseteq A$ of square zero, the surjectivity of $G(A)\longrightarrow G(A/I)$ is clear from the definition of $G$, hence $G/\Spec(\Oh_k)$ is 
smooth. By Proposition \ref{flat1}, {\sl ii)}$\Rightarrow${\sl i)}, $\Nr:G\longrightarrow \G_m$ is (faithfully) flat, hence so is its base change $G'/\Spec(\Oh_k)$.\\
We shall now show that all geometric fibers of $G$ (resp. $G'$) are
smooth and connected of dimension $d^2$ (resp. $d^2-1$). This will also imply that $G'/\Spec(\Oh_k)$ is smooth by the fiber-wise criterion \cite[Th\'eor\`eme 17.5.1]{ega44} and thus conclude the proof.\\
Geometric fibers of $G$ (resp. $G'$) in
characteristic zero are isomorphic to $\Gl_d$ (resp. $\Sl_d$). Let $0\neq\p\subseteq\Oh_k$ be a prime, $\kappa:=\Oh_k/\p$ and $\overline{\kappa}$ be an algebraic closure of
$\kappa$. We have $D_{\p}\simeq\M_n(\Df)$ for a central skew-field $\Df$ over $k_{\p}$. Denoting by $r$ the reduced dimension of $\Df$, we have $d=nr$. Since
$\Oh\otimes_{\Oh_k}\Oh_{k,\p}\subseteq D_{\p}$ is a maximal order \cite[Corollary 11.2]{reiner}, we have $\Oh\otimes_{\Oh_k}\Oh_{k,\p}\simeq\M_n(\Oh_{\Df})$ as $\Oh_{k,\p}$- algebras by \cite[Theorem 17.3]{reiner} where
$\Oh_{\Df}\subseteq\Df$ is the unique maximal order \cite[Theorem 12.8]{reiner}.\\
Let $\Pi\in\Oh_{\Df}$ and $\pi\in\Oh_{k,\p}$ be uniformizers. Then $\overline{A}:=
(\Oh_{\Df}/\pi\Oh_{\Df})\otimes_{\kappa}\overline{\kappa}$ is a $\overline{\kappa}$- 
algebra with radical $\overline{\Rf}=(\Pi\Oh_{\Df}/\pi\Oh_{\Df})\otimes_{\kappa}\overline{\kappa}$ and maximal semi-simple quotient $\overline{A}/\overline{\Rf}\simeq\overline{\kappa}^r$. Since $G_{\overline{\kappa}}=\Gl_n(\overline{A})$, we 
have an extension

\[ \xymatrix{ 1 \ar[r] & U \ar[r] & G_{\overline{\kappa}} \ar[r] & (\Gl_{n,\overline{\kappa}})^r \ar[r] & 1 } \]

where $U$ is a unipotent group of dimension $n^2(r-1)r$ (recall that 
$\pi\Oh_{\Df}=\Pi^r\Oh_{\Df}$ and $(\Pi^i\Oh_{\Df}/\Pi^{i+1}\Oh_{\Df})\otimes_{\kappa}\overline{\kappa}\simeq\overline{\kappa}^r$). So $G_{\overline{\kappa}}$ is connected and smooth of dimension $n^2r+n^2(r-1)r=d^2$.\\
Since the reduced norm $\Nr_{\overline{\kappa}}:G_{\overline{\kappa}}\longrightarrow\G_{m,\overline{\kappa}}$ is trivial on $U$ it factors over some $\alpha:(\Gl_{n,\overline{\kappa}})^r\longrightarrow\G_{m,\overline{\kappa}}$. We have $\alpha(g_1,\ldots,g_r)=\prod_{i=1}^r\det(g_i)$ as an immediate consequence of \cite[Lemma 3.8]{kleinert}. This exhibits $G_{\overline{\kappa}}'$ as an extension of $V:=\ker(\alpha)$ by $U$. We can factor $\alpha=\beta\circ\gamma$ with $\gamma:(\Gl_{n,\overline{\kappa}})^r\longrightarrow \G_{m,\overline{\kappa}}^r$, $\gamma(g_1,\ldots, g_r):=(\det(g_i))_i$ and $\beta:\G_{m,\overline{\kappa}}^r\longrightarrow\G_{m,\overline{\kappa}}$, $\beta(x_1,\ldots, x_r):=x_1\ldots x_r$ and thus obtain

\[ \xymatrix{ & & & 1 \ar[d] & \\ 1 \ar[r] & \Sl_{n,\overline{\kappa}}^r \ar[r]\ar[d]^{\id} & V \ar[r]\ar[d]  & T\ar[r] \ar[d] & 1\\
1 \ar[r] & \Sl_{n,\overline{\kappa}}^r \ar[r]  &  \Gl_{n,\overline{\kappa}}^r \ar[r]^{\gamma} \ar[dr]^{\alpha}  & \G_{m,\overline{\kappa}}^r \ar[r] \ar[d]^{\beta}  & 1 \\  & & & \G_{m,\overline{\kappa}} \ar[d] &  \\
 & & & 1 & } \]

with $T:=\ker(\beta)$. Looking at character groups for example, one sees that 
$T\simeq\G_{m,\overline{\kappa}}^{r-1}$ and hence $V$ is connected and smooth 
of dimension $\dim(T)+\dim(\Sl_n^r)=n^2r-1$ and $G'_{\overline{\kappa}}$ is 
connected and smooth of dimension $\dim(V)+\dim(U)=d^2-1$.

\end{proofof}

\begin{remark} The maximal locus inside $\Spec(\Oh_k)$ over which $G$ (resp. $G'$) is reductive (resp. semi-simple) is obtained by inverting the discriminant of $D$, i.e. by removing all $\p\in\Spec(\Oh_k)$ such that $\inv_{\p}(D)\neq 0$.
\end{remark}

\subsubsection{Type $^2A_{d-1}$}\label{geo2ad-1}\hfill

Let $D$ be a finite-dimensional skew-field of reduced dimension $d$ over $\Q$ carrying a positive involution 
of the second kind $*$, i.e. for all $x\in D^*$ we have $\tr^D_{\Q}(^* xx)>0$ (positivity)
and $*$ restricted to the center $L$ of $D$ is non-trivial. Then $L$ is a 
CM-field with $k:=\{ x\in L\, | \, x=^*x \}\subseteq L$ as its maximal real
subfield \cite[page 194]{mumford}. Note that $*$ is $k$-linear. We assume that 
$\Oh\subseteq D$ is a maximal order which is invariant under $*$. The existence of such an order is claimed without proof in \cite[7.1.1]{hida}. Then $\Oh\cap
L=\Oh_L$ and $\Oh\cap k=\Oh_k$ are the rings of integers of $L$ and $k$.
We consider the affine group-schemes of finite type $\U$ and $\T$ over $\Spec(\Oh_k)$ whose groups of points are given
for every $\Oh_k$-algebra $R$ by

\[ \U(R)=\{ g\in(\Oh\otimes_{\Oh_k} R)^*\, | \, ^*gg=1 \}\mbox{   and} \]
\[ \T(R)=\{ g\in(\Oh_L \otimes_{\Oh_k} R)^*\, | \, \Nr^L_k(g)=1 \}. \]

There is a homomorphism $\Nr:\U\longrightarrow \T$ over $\Spec(\Oh_k)$ given on points by the reduced norm of $D$ and we put $\SU:=\ker(\Nr)$.
Over $\Spec(k)$ we have an exact sequence

\[ \xymatrix{ 1\ar[r] & \SU_1(D,1)=\SU\times_{\Spec(\Oh_k)}\Spec(k)  \ar[r] & \U_1(D,1)=\U\times_{\Spec(\Oh_k)}\Spec(k)  \ar[r]^-{\Nr_k} & \Res^L_k(\G_{m,L})^{(1)} \ar[r]& 1,} \]

where "1" denotes the standard rank one Hermitian form on $D$ and \[ \Res^L_k(\G_{m,L})^{(1)}:=\ker(\Res^L_k(\G_{m,L})\stackrel{\Nr_k^L}{\longrightarrow}\G_{m,k})\, \] is a one-dimensional anisotropic torus over $k$; c.f. \cite[Section 2.3]{pr}
for notation and general background on unitary groups.\\
We first study the integral model $\T/\Spec(\Oh_k)$ of $\Res^L_k(\G_{m,L})^{(1)}$. We define the open subscheme ${\mathcal U}\subseteq\Spec(\Oh_k)$ by

\[ {\mathcal U}:=\Spec(\Oh_k)-\{ 0\neq\p\subseteq\Oh_k\mbox{ is a prime of residue characteristic } 2\mbox{ and ramified in } L/k\}. \]

The following result makes \cite[Proposition 5.2]{sansuc} slightly more precise in the present special case.

\begin{prop}\label{geoofT} $\T/\Spec(\Oh_k)$ is an affine flat group-scheme such that $\T_k\simeq \Res^L_k(\G_{m,L})^{(1)}$. For a prime $0\neq\p\subseteq\Oh_k$ we have
\[ \T_{\kappa(\p)}\simeq\left\{ \begin{array}{ccc} \G_{m,\kappa(\p)} & , & \mbox{ if }\p\mbox{ splits in }L/k.\\
\Res^{\kappa(\p)^{(2)}}_{\kappa(\p)}(\G_{m,\kappa(\p)^{(2)}})^{(1)} & , & \mbox{ if }\p\mbox{ is inert in }L/k.\\
\G_{a,\kappa(\p)}\times\mu_{2,\kappa(\p)} & , & \mbox{ if }\p\mbox{ is ramified in }L/k.\end{array}\right. \]
In particular, the maximal locus inside $\Spec(\Oh_k)$ over which $\T$ is smooth equals ${\mathcal U}$.\\
Here, $\kappa(\p):=\Oh_k/\p$ and $\kappa(\p)^{(2)}$ is the 
unique quadratic extension of $\kappa(\p)$.
\end{prop}

\begin{proof} We know that $\Res^{\Oh_L}_{\Oh_k}(\G_{m,\Oh_L})/\Spec(\Oh_k)$ is an affine and smooth group-scheme from \cite[7.6, proof of Theorem 4 and Proposition 5,h)]{blr}. There is an obvious subgroup $i:\G_{m,\Oh_k}\hookrightarrow \Res^{\Oh_L}_{\Oh_k}(\G_{m,\Oh_L})$ such that 
$\Nr^L_k\circ i$ is multiplication by $2$, hence $\Nr^L_k: \Res^{\Oh_L}_{\Oh_k}(\G_{m,\Oh_L})\longrightarrow \G_{m,\Oh_k}$ is an $fppf$-epimorphism and the first assertion follows from Proposition
\ref{flat1} since by definition $\T=\ker(\Nr_k^L)$. Since restriction commutes with base change, for every $\Oh_k$-algebra $R$ we have 

\[ \T_R=\ker(\Res^{\Oh_L\otimes_{\Oh_k} R}_R(\G_{m,R})\longrightarrow \G_{m,R})\]

which makes the assertions concerning the generic fiber and the fibers over split and inert primes obvious.\\
For ramified $\p$ we have $\Oh_L\otimes_{\Oh_k}\kappa(\p)\simeq \Oh_{L,\qf}/\qf^2\Oh_{L,\qf}$ for $\qf$ the unique prime of $\Oh_L$ lying above $\p$. There exists $\alpha\in \Oh_{L,\qf}$ with $\Oh_{L,\qf}=\Oh_{k,\p}[\alpha]$ and $\alpha$ satisfies an Eisenstein polynomial $x^2-ax+b\in\Oh_{k,\p}[x]$. Since $a\in\p\Oh_{k,\p}\subseteq \qf^2\Oh_{L,\qf}$, the non trivial automorphism $\sigma$ of $\Oh_{L,\qf}$ over $\Oh_{k,\p}$ satisfies

\[ \sigma(\alpha)=a-\alpha\equiv -\alpha\mbox{ in }\Oh_{L,\qf}/\qf^2\Oh_{L,\qf}.\]

As $\Oh_{L,\qf}/\qf^2\Oh_{L,\qf}\simeq\kappa(\p)[\alpha]/(\alpha^2)$ we conclude that for every $\kappa(\p)$-algebra $R$

\[ \T_{\kappa(\p)}(R)\simeq\{ x+y\alpha\in(R[\alpha]/(\alpha^2))^*\, | \, 1=(x+y\alpha)\sigma(x+y\alpha)=(x+y\alpha)(x-y\alpha)=x^2\} \]

and we have an exact sequence

\[ \xymatrix{ 1 \ar[r] & \G_{a,\kappa(\p)}(R) \ar[r] & \T_{\kappa(\p)}(R) \ar[r] & \mu_{2,\kappa(\p)}(R) \ar[r] & 1.\\
 & & x+y\alpha \ar@{|->}[r] & x & } \]

which is split by $x\mapsto x+0\alpha$.

\end{proof}

We will need the following.

\begin{prop}\label{switch}
Let $k$ be a commutative ring, $B_1$ and $B_2$ $k$-algebras and $\tau$ an involution on  $B:=B_1\times B_2$ such that $\tau(x,y)=(y,x)$ for all $x,y\in k$. Then there is an isomorphism of $k$-algebras with involution

\[ (B,\tau)\simeq (B_1\times B_1^{opp}, (x,y)\mapsto (y,x)).\]
\end{prop}

\begin{proof} The proof of \cite[Proposition 2.14]{bookinv} carries over without any change.
\end{proof}

Now let $C$ be the set of non-zero primes $\p\in {\mathcal U}$ such that $\U_{\overline{\kappa(\p)}}$ is an extension of a {\em symplectic} group. We will see during the proof of Theorem \ref{geoofUandSU} that this set only contains primes which are ramified in $L/k$. Let $\T'\subseteq \T$ be the open subscheme obtained by removing from $\T$ the non-identity component of $\T_{\kappa(\p)}$ for all $\p\in C$, c.f. Proposition \ref{geoofT}. Clearly, $\T'\subseteq\T$ is a subgroup-scheme.

\begin{theorem}\label{geoofUandSU}
i) The morphism $\Nr_{{\mathcal U}}:\U_{{\mathcal U}}\longrightarrow\T_{{\mathcal U}}$ factors through $\T'_{{\mathcal U}}\subseteq\T_{{\mathcal U}}$ and the resulting sequence of $fppf$-sheaves on ${\mathcal U}$
\[ 1\longrightarrow \SU_{{\mathcal U}}\longrightarrow\U_{{\mathcal U}}\stackrel{\Nr_{{\mathcal U}}}{\longrightarrow}\T'_{{\mathcal U}}\longrightarrow 1\] is exact.\\
ii) The group  $\U/{\mathcal U}$ (resp. $\SU/{\mathcal U}$) is smooth (resp. smooth with connected fibers).\\

\end{theorem}

\begin{remark}
It is easy to give examples of the situation in Theorem \ref{geoofUandSU} in which $C\neq\emptyset$, i.e. $\T'\neq \T$. From case $3'.2)$ in the proof of Theorem \ref{geoofUandSU} it will however be clear that $C=\emptyset$ if $d$ is odd.
\end{remark}

\begin{proof} Fix $\p\in {\mathcal U}$. We will study the local groups 
\[ \U_{\p}:=\U\times_{\Spec(\Oh_k)}\Spec(\Oh_{k,\p})\mbox{ and }\]
\[ \SU_{\p}:=\SU\times_{\Spec(\Oh_k)}\Spec(\Oh_{k,\p})\]
over $\Spec(\Oh_{k,\p})$. For this, we need to understand the $\Oh_{k,\p}$-algebra with involution 

\[ \Oh_{\p}:=\Oh\otimes_{\Oh_k}\Oh_{k,\p}.\]

We distinguish three cases:\\

1) $\p$ is inert in $L/k$:\\
For the unique prime $\qf\subseteq\Oh_L$ lying above $\p$
we have $\inv_{\qf}(D)=0$ \cite[(B) on page 199]{mumford} and

\[ \Oh_{\p}\subseteq\Oh_{\p}\otimes_{\Oh_{k,\p}} k_{\p}\simeq D\otimes_L L_{\qf}\simeq \M_d(L_{\qf}) \]

is a maximal order, hence $\Oh_{\p}\simeq 
\M_d(\Oh_{L,\qf})$ as $\Oh_{k,\p}$-algebras. Denote by $\sigma$ the non-trivial 
automorphism of $\Oh_{L,\qf}$ over $\Oh_{k,\p}$. 
We obtain for every $\Oh_{L,\qf}$-algebra $R$:

\[ (\U_{\p}\times_{\Spec(\Oh_{k,\p})} \Spec(\Oh_{L,\qf}))(R)=\U_{\p}(R)\simeq \{ g\in\M_d(\Oh_{L,\qf}\otimes_{\Oh_{k,\p}} R)\, | \, ^*gg=1\}.\]

Since $\p$ is unramified we have $\Oh_{L,\qf}\otimes_{\Oh_{k,\p}}\Oh_{L,\qf}\simeq \Oh_{L,\qf}\times\Oh_{L,\qf}$ and
under this isomorphism $(\sigma\otimes 1)$ switches the factors.
Since by Proposition \ref{switch} we have an isomorphism of $\Oh_{L,\qf}$-algebras with involution

\[ (\M_d(\Oh_{L,\qf})\otimes_{\Oh_{k,\p}}\Oh_{L,\qf},*)\simeq (\M_d(\Oh_{L,\qf})\times \M_d(\Oh_{L,\qf})^{opp}, (x,y)\mapsto (^ty,^tx)),\]

where $^t$ denotes the transpose of a matrix, we find that

\[ \U_{\p}(R)\simeq\{ (x,y)\in \M_d(R\times R)\simeq \M_d(R)\times\M_d(R)\, | \, (^ty,^tx)(x,y)=1\}\simeq\Gl_d(R).\]

We have thus shown that $\U_{\p}\times_{\Spec(\Oh_{k,\p})}\Spec(\Oh_{L,\qf})\simeq\Gl_{d,\Oh_{L,\qf}}$. By decent, we conclude that 
$\U_{\p}/\Spec(\Oh_{k,\p})$ is smooth. Furthermore, the 
special fiber $\U_{\kappa(\p)}:=\U_{\p}\times_{\Spec(\Oh_{k,\p})}\Spec(\kappa(\p))$ is a $\overline{\kappa(\p)}/\kappa(\p)$- form of $\Gl_{d,\kappa(\p)}$ and since $\HH^1(\Spec(\kappa(\p)),\PGL_d)=1$ we have $\U_{\kappa(\p)}\simeq\Gl_{d,\kappa(\p)}$.\\

2) $\p$ splits in $L/k$:\\
In this case we have \[ \Oh_{\p}\simeq \Oh_{\qf}\times\Oh_{\qf'}\]
where $\qf$ and $\qf'$ are the primes of $\Oh_L$ lying above $\p$
and $\Oh_{\qf}:=\Oh\otimes_{\Oh_L}\Oh_{L,\qf}$ and similarly for $\qf'$. By Proposition \ref{switch} we have an
isomorphism of $\Oh_{k,\p}$-algebras with involution 

\[ (\Oh_{\p},*)\simeq(\Oh_{\qf}\times\Oh_{\qf}^{opp},(x,y)\mapsto (y,x))\]

and thus $\U_{\p}\simeq\Gl_1(\Oh_{\qf})$ and this group is trivially smooth over $\Spec(\Oh_{k,\p})$.\\

3) $\p$ is ramified in $L/k$: as in case 1) we have

\[ \Oh_{\p}\simeq\M_d(\Oh_{L,\qf})\]

as $\Oh_{k,\p}$-algebras with $\qf\subseteq\Oh_L$ the unique prime
lying above $\p$. We check the
smoothness of $\U_{\p}/\Spec(\Oh_{k,\p})$ using the lifting criterion: let $A$ be an Artinian $\Oh_{k,\p}$-algebra
and $I\subseteq A$ an ideal with $I^2=0$. Given

\[ x\in\U_{\p}(A/I)=\{ g\in(\Oh_{\p}\otimes_{\Oh_{k,\p}} A/I)^*\, | \,
^*gg=1\} \]

there is some $y\in \Oh_{\p}\otimes_{\Oh_{k,\p}} A$ lifting $x$ and we
have

\[ ^*yy=1+z\mbox{ for some } z\in\Oh_{\p}\otimes_{\Oh_{k,\p}}
I\subseteq\Oh_{\p}\otimes_{\Oh_{k,\p}} A  \]

with $^*z=z$ because $^*(^*yy)=^*yy$. As $\p\in {\mathcal U}$ we have $2\in
A^*$ and can define

\[ y':=y(1-\frac{1}{2}z)\in\Oh_{\p}\otimes_{\Oh_{k,\p}} A\]

which still lifts $x$ and satisfies

\[ ^*y'y'=(1-\frac{1}{2}\, ^*z)^*yy(1-\frac{1}{2}z)=(1-\frac{1}{2}\,^*z)(1+z)(1-\frac{1}{2}z)\stackrel{(I^2=0,^*z=z)}{=}1. \]

Hence we have found $y'\in\U_{\p}(A)$ lifting $x$.\\

At this point we have established that $\U/{\mathcal U}$
is smooth and now proceed to study $\SU/{\mathcal U}$.
We first consider the geometric fibers showing in particular that they
are all connected and smooth.\\
Let $\Spec(\Omega)\longrightarrow {\mathcal U}$ be a geometric point. If the
characteristic of $\Omega$ is zero, we have
$\SU_{\Omega}\simeq\Sl_{d,\Omega}$, hence we can assume that
$\Omega=\overline{\kappa(\p)}$ for some $\p\in {\mathcal U}$. We again have
to distinguish three cases as above:\\

1') $\p$ is inert in $L/k$: from 1) above we have
$\U_{\Omega}\simeq\Gl_{d,\Omega}$ and from Proposition \ref{geoofT} we
know that $\T_{\Omega}\simeq\G_{m,\Omega}$. Under these isomorphisms,
$\Nr_{\Omega}$ is identified with
the determinant, hence $\SU_{\Omega}\simeq\Sl_{d,\Omega}$. \\

2') $\p$ splits in $L/k$: from 2) above and Proposition \ref{geoofT} we
know that $\U_{\Omega}\simeq\Gl_1(\overline{A})$ and
$\T_{\Omega}\simeq\G_{m,\Omega}$ for the $\Omega$-algebra $\overline{A}:=\Oh_{\qf}\otimes_{\Oh_{k,\p}}\Omega$ where $\qf\subseteq\Oh_L$ is a prime lying above $\p$. This is the
situation studied in Theorem \ref{smooth} from which we read off that
$\SU_{\Omega}/\Spec(\Omega)$ is connected and smooth (and, in fact, also the
dimensions of the semi-simple, toric and unipotent parts of $\SU_{\Omega}$ in terms of the order of $\inv_{\qf}(D)$).\\

3') $\p$ is ramified in $L/k$: recall from 3) above that
$\Oh_{\p}\simeq\M_d(\Oh_{L,\qf})$. We have to study the involution induced
by $*$ on

\[ \Oh_{\p}\otimes_{\Oh_{k,\p}}\Omega\simeq\M_d(\Omega(\epsilon)).\]

Recall that $\Omega(\epsilon):=\Omega[\epsilon]/(\epsilon^2)$ and that
the involution $*$ on
$\M_d(\Omega(\epsilon))$ satisfies $^*\epsilon=-\epsilon$ as established
during the proof of Proposition \ref{geoofT}. Denoting by
$\sigma\in\Aut_{\Omega-alg}(\Omega(\epsilon))$ the
element determined by $\sigma(\epsilon)=-\epsilon$ and by $+$ the
involution on $\M_d(\Omega(\epsilon))$ defined by

\[ ^+(x_{i,j}):=(\sigma(x_{j,i})),\]

the Theorem of Skolem-Noether \cite[Chapter IV, Proposition 1.4]{milneetlae} shows that there exists a
$g\in\Gl_d(\Omega(\epsilon))$ such that
\begin{equation}\label{inv2}
^*x=g ^+x g^{-1}\mbox{ for all }x\in\M_d(\Omega(\epsilon)).
\end{equation}

From $^{**}x=x$ one sees that 

\begin{equation}\label{abc}
g=\alpha ^+g
\end{equation}

for some $\alpha\in\Omega(\epsilon)^*$. This gives $^+g=\sigma(\alpha) g$ and by multiplying we obtain $g^+g=\alpha\sigma(\alpha) ^+gg$ which using (\ref{abc}) implies that $\alpha\sigma(\alpha)=1$.\\
Writing $\alpha=x+y\epsilon$ with $x,y\in\Omega$ we get

\[ 1=\alpha\sigma(\alpha)=x^2-y^2\epsilon^2=x^2\]

hence

\begin{equation}\label{alpha}
\alpha=\pm  1+y\epsilon\mbox{ for some }y\in\Omega.
\end{equation}

Replacing $g$ by $\beta g$ for
$\beta:=1\mp\frac{1}{2}y\epsilon\in\Omega(\epsilon)^*$ (the
sign opposite to the one occurring in (\ref{alpha}))
does not affect (\ref{inv2}) and replaces $\alpha$ by

\[ \alpha\beta\sigma(\beta)^{-1}\stackrel{(\ref{alpha})}{=}(\pm
1+y\epsilon)(1\mp\frac{1}{2}y\epsilon)(1\pm\frac{1}{2}y\epsilon)^{-1}=(\pm
1+y\epsilon)(1\mp\frac{1}{2}y\epsilon)^2=\]
\[ (\pm 1+y\epsilon)(1\mp y\epsilon)=\pm 1.\]

Hence we can assume that

\begin{equation}\label{alpha2}
\alpha=\pm 1.
\end{equation}

To further simplify the involution, note that for every
$h\in\Gl_d(\Omega(\epsilon))$, $(\M_d(\Omega(\epsilon)),*)$ is
isomorphic, via conjugation with $h$, to $(\M_d(\Omega(\epsilon)),\tau)$
with

\[ ^{\tau}x=h^*(h^{-1}xh)h^{-1}=hg^+(h^{-1}xh)g^{-1}h^{-1}= hg^+h
^+x (hg^+h)^{-1},\]

i.e.

\begin{equation}\label{replace}
\mbox{ For every }h\in\Gl_d(\Omega(\epsilon))\mbox{ we can replace }g\mbox{ in }(\ref{inv2})\mbox{ by }hg^+h.
\end{equation}

We now distinguish two cases according to (\ref{alpha2}):\\

3'.1) Assume that $\alpha=1$. Writing $g=A+B\epsilon$ with
$A\in\Gl_d(\Omega),
B\in\M_d(\Omega)$ we have

\[ A+B\epsilon=g\stackrel{((\ref{abc}),\alpha=1)}{=}\, ^+g=^tA-
^tB\epsilon,\]

hence $A=^tA,B=-^tB$ and there exists some $S\in\Gl_d(\Omega)$ with
$A=S^tS$. Using (\ref{replace}) with $h=S^{-1}$ we can
replace $g$ by

\[ hg^+h=S^{-1}(A+B\epsilon)^tS^{-1}=1+S^{-1}B^tS^{-1}\epsilon.\]

Put $T:=S^{-1}B^tS^{-1}$ and note that $B=-^tB$ implies that $^tT=-T$.
Using
(\ref{replace}) again with $h=1-\frac{1}{2}T\epsilon$ replaces $g$ by

\[
hg^+h=(1-\frac{1}{2}T\epsilon)(1+T\epsilon)(1+\frac{1}{2}\, ^tT\epsilon)=1,\]

i.e. we can assume that $^*x=^+x$ for all $x\in\M_d(\Omega(\epsilon))$.
For every $\Omega$-algebra $R$ we thus obtain

\[ \U_{\Omega}( R ) \simeq\{ x=X+Y\epsilon\in\M_d(R(\epsilon))\, | \,
1=^+xx=(^tX-^tY\epsilon)(X+Y\epsilon)= \]
\[ ^tXX+(-^t YX+^tXY)\epsilon \}\]

and hence an exact sequence

\[ \xymatrix{ 1\ar[r] & F( R):=\{ 1+Y\epsilon\in\M_d(\Omega(\epsilon))\,
| \, Y=^tY \} \ar[r] & \U_{\Omega}( R) \ar[r] & \O_d( R)\ar[r] & 1\\ & &
X+Y\epsilon\ar@{|->}[r] & X }, \]

$\O_d$ denoting the orthogonal group, which is split by $X\mapsto X+0\cdot\epsilon$. Evidently,
$F\simeq\G_{a,\Omega}^{\frac{1}{2}d(d+1)}$. We have the following
diagram with exact rows and columns:

\begin{equation}\label{diagram1}
\xymatrix{ & 1\ar[d] &1 \ar[d] & 1 \ar[d]  & \\ 1 \ar[r] & F'\ar[r]
\ar[d] & \SU_{\Omega} \ar[r]^{\alpha}  \ar[d] & \SO_{d,\Omega} \ar[r]
\ar[d] &1 \\ 1 \ar[r] & F \ar[r]^{\iota} \ar[d]^{\tr} & \U_{\Omega}
\ar[r] \ar[d]^{\Nr_{\Omega}} & \O_{d,\Omega} \ar[r] \ar[d]^{\det} & 1
\\
1 \ar[r] & \G_{a,\Omega}\ar[d] \ar[r] & \T_{\Omega} \ar[d] \ar[r]^{\pi} & \mu_{2,\Omega}
\ar[r] \ar[d] & 1\\
 & 1 & 1 & 1 & }
\end{equation}

This is obtained as follows: the lower row is taken from Proposition
\ref{geoofT}. The reduced norm induces the determinant on
$\M_d(\Omega(\epsilon))$. This shows that $\pi\Nr_{\Omega}\iota$ is
trivial and $\Nr_{\Omega}\iota$ factors through some
$F\longrightarrow\G_a$. As $\det(1+Y\epsilon)=1+\tr(Y)\epsilon$, the
map $F\longrightarrow \G_a$ is in fact the trace and $F':=\ker(\tr)$.
This also shows that the map $\O_{d,\Omega}\longrightarrow\mu_2$
induced by $\Nr_{\Omega}$ is the determinant which is visibly $fppf$-surjective, in fact, it is surjective as a morphism of pre-sheaves as is the trace $\tr$. It is also clear that
$F'\simeq\G_a^{\frac{1}{2}d(d+1)-1}$. Now the $fppf$-surjectivity of $\alpha$
and $\Nr_{\Omega}$ follows from a $5$-lemma argument (which does not use commutativity).
In particular, $\SU_{\Omega}$ is connected and smooth.\\

3'.2) Assume that $\alpha=-1$. Writing $g=A+B\epsilon$ with
$A\in\Gl_d(\Omega), B\in\M_d(\Omega)$ we have

\[
A+B\epsilon=g\stackrel{((\ref{abc}),\alpha=-1)}{=}-^+g=-^tA+^tB\epsilon,\]

i.e. $A=-^tA$ and $B=^tB$. The conditions on $A$ force $d$ to be even,
say $d=2m$. Let
$J:=\left(\begin{array}{cc} 0 & 1_m\\-1_m &
0\end{array}\right)\in\Gl_d(\Omega)$ be the standard alternating
matrix. Then there exist a $S\in\Gl_d(\Omega)$ such that $SA^tS=J$. Using (\ref{replace}) with $h=S$ we can replace $g$ by

\[ hg^+h=S(A+B\epsilon)^tS=J+SB^tS\epsilon.\]

Put $T:=SB^tS$ and note that $B=^tB$ implies that $T=^tT$. Using
(\ref{replace}) again with $h=1+\frac{1}{2}TJ\epsilon$ replaces $g$ by

\[
hg^+h=(1+\frac{1}{2}TJ\epsilon)(J+T\epsilon)(1-\frac{1}{2}\,^tJ^tT\epsilon)=J+(\frac{1}{2}TJ^2+T-\frac{1}{2}J^tJ^tT)\epsilon=\]
\[ J+(\frac{1}{2}T(-1)+T-\frac{1}{2}\, ^tT)\stackrel{(^tT=T)}{=}J,\]

i.e. we can assume that $g=J$.\\
For every $\Omega$-algebra $R$ we thus obtain, using $^tJ=J^{-1}$,

\[ \U_{\Omega}( R)\simeq \{ x=X+Y\epsilon\in\M_d(\Omega(\epsilon))\, | \,
1=^*xx=J^+(X+Y\epsilon)^tJ(X+Y\epsilon)=\]
\[ J(^tX-^tY\epsilon)^t
J(X+Y\epsilon)=J^tX^tJX+(-J^tY^tJX+J^tX^tJY)\epsilon\ \}.\]

Note that $1=J^tX^tJX$ if and only if $^tXJX=J$, hence we get an exact
sequence

\[ \xymatrix{ 1 \ar[r] & F( R):=\{
1+Y\epsilon\in\M_d(\Omega(\epsilon))\, | \, Y=J^tY^tJ\} \ar[r] &
\U_{\Omega}( R) \ar[r] & \Sp_{2m}( R)\ar[r] & 1 \\
 & & X+Y\epsilon \ar@{|->}[r] & X & }, \]

where $\Sp_{2m}$ denotes the symplectic group, which is split by $X\mapsto X+0\epsilon$. Writing $Y=\left(
\begin{array}{cc} a & b\\c & d\end{array}\right)$ with $a,b,c,d\in\M_m(\Omega)$ one obtains

\[ F( R)\simeq\left\{ \left(
\begin{array}{cc} a & b \\ c & d \end{array}\right)\in\M_{2m}( R)\, | \, ^ta=d,
^tb=-b, ^tc=-c\right\},\]

hence $F\simeq\G_a^{2m^2-m}$. The analogue of diagram
(\ref{diagram1}) in this case reads as follows.

\begin{equation}\label{diagram2}
\xymatrix{ & 1 \ar[d] & 1\ar[d] & 1\ar[d] & \\ 1 \ar[r] & F'\ar[r]\ar[d] &
\SU_{\Omega}\ar[d]\ar[r]^{\alpha} & \Sp_{2m}\ar[r]\ar[d] & 1\\ 1 \ar[r]  & F
\ar[r]\ar[d]^{\tr} & \U_{\Omega} \ar[r]\ar[d]^{\Nr_{\Omega}} \ar@{.>}[dl] & \Sp_{2m}
\ar[r]\ar[d]^{\det(\equiv 1)} & 1\\
1\ar[r] & \G_{a,\Omega} \ar[d]\ar[r] & \T_{\Omega} \ar[r] & \mu_{2,\Omega} \ar[r] & 1 \\ & 1 & & &}
\end{equation}

Note that since the determinant of every symplectic matrix equals $1$, $\Nr_{\Omega}$ factors as indicated in diagram (\ref{diagram2}). In particular, $\Nr_{\Omega}:\U_{\Omega}\longrightarrow\T_{\Omega}$ is not an $fppf$-epimorphism but has image the connected component $\G_{a,\Omega}\simeq\T^0_{\Omega}\subseteq \T_{\Omega}$.
Since $F'\simeq \G_a^{2m^2-m-1}$ we conclude that $\SU_{\Omega}$ is connected and smooth.\\

We now establish the exactness of the sequence 

\[ 1\longrightarrow \SU_{{\mathcal U}} \longrightarrow \U_{{\mathcal U}} \stackrel{\Nr_{{\mathcal U}}}{\longrightarrow} \T'_{{\mathcal U}} \longrightarrow 1\]

over ${\mathcal U}$. Since $\T'\subseteq \T$ is an open subscheme, the fact that $\Nr_{{\mathcal U}}$ factors through $\T'_{{\mathcal U}}$ can be checked on fibers where it follows from what has been shown above: since the determinant is trivial on $\Sp_{2m}$ in diagram (\ref{diagram2}), $\Nr_{\Omega}$ factors through $\G_{a,\Omega}=\T'_{\Omega}\subseteq \T_{\Omega}$. We now need to see that the resulting morphism $\Nr_{{\mathcal U}}: \U_{{\mathcal U}}\longrightarrow \T'_{{\mathcal U}}$ is an $fppf$-epimorphism and we will show that it is in fact faithfully flat: by Proposition \ref{flat1}, it is enough to see that $\U_{\Omega}\longrightarrow \T'_{\Omega}$ is an 
$fppf$-epimorphism for every geometric point $\Spec(\Omega)\longrightarrow {\mathcal U}$ which follows again from what has been proved above. The flatness of $\U_{{\mathcal U}}\longrightarrow\T'_{{\mathcal U}}$ implies that $\SU/{\mathcal U}$ is flat, hence smooth by the
fiber-wise criterion.\\
The proof of Theorem \ref{geoofUandSU} is now complete.

\end{proof}

\section{Approximation of local units}\label{approxlocal}

In this section we study the problem of $p$-adically approximating
local units of a maximal order (with involution) by global (unitary) units of bounded denominators.
Using the results of subsections \ref{torsor} and \ref{geoad-1} (\ref{geo2ad-1})
this problem will be reduced in subsection \ref{ad-1} (\ref{2ad-1}) to an approximation problem for tori which will be solved in subsection \ref{cad-1} (solved in a special case in subsection \ref{c2ad-1}).

\subsection{Type $A_{d-1}$}\label{ad-1}

In this subsection we consider the problem of $p$-adically approximating local units
of a maximal order $\Oh\subseteq D$ where $D$ is a finite dimensional skew-field
over $\Q$. We denote by $k$ the center of $D$ and by $d$ its reduced dimension. 
We fix a prime $0\neq\p\subseteq\Oh_k$ at which we wish to approximate. 
There is a unique prime $\P\subseteq\Oh$ lying above $\p$ \cite[VI, \S 12, Satz 1]{deu} and we denote by
$\Oh_{\P}$ the $\P$-adic completion of $\Oh$, c.f. \cite[Kapitel VI, \S 11]{deu}.\\
To describe the denominators we allow the approximating global units to have, we fix a finite set of places $\Srm$
of $k$ such that
\[ \p\not\in\Srm \mbox{ and there exists a place } v_0\in\Srm\mbox{ such that }D_{v_0}\mbox{ is not a skew-field.} \]
We write $\Srm^{fin}$ for the set of finite places contained in $\Srm$ and consider the ring $\OSfin$ of $\Srm^{\mathrm fin}$-integers

\[ \Oh_k\subseteq\OSfin:=\{ x\in k\st v(x)\ge 0\mbox{ for all finite }v\not\in\Srm \}\subseteq k.\]

Since $\p\not\in\Srm$ we have $\OSfin\subseteq\Oh_{k,\p}$. We define 

\begin{equation}\label{Xdef}
 X:=\{ x\in\OSfin^*\st v\mbox{ infinite and }\inv_v(D)=\frac{1}{2}\mbox{ imply }v(x)>0\}\subseteq\Oh_{k,\p}^*\mbox{  and }
\end{equation}

\[ (\Oh\otimes_{\Oh_k}\OSfin)^*\subseteq\Oh_{\P}^*.\]

Recall that $\Nr$ denote the reduced norm of $D$.

\begin{theorem}\label{approx} The closure of $(\Oh\otimes_{\Oh_k}\OSfin)^*$ inside
$\Oh_{\P}^*$ equals \[ \{ x\in\Oh_{\P}^*\st \Nr_{\p}(x)\in\Oh_{k,\p}^*\mbox{ lies in the closure of }X \}. \]
\end{theorem}
\vspace{0.2cm}
\begin{example}\label{tralala}
1) For $k=\Q$ and $D$ a definite quaternion algebra, i.e. $d=2$ and 
$\inv_v(D)=\frac{1}{2}$ for the unique infinite place $v$ of $\Q$, we can choose
$\Srm=\{ l \}$ for any prime $l\neq p$ at which $D$ splits, i.e. $\inv_l(D)=0$.
Then $\OSfin^*=\{\pm 1 \}\times l^{\Z}$ and $X=l^{\Z}\subseteq\Oh_{k,\p}^*=\Z_p^*$.
For $p\neq 2$ we can choose $l$ as above such that in addition $X\subseteq\Z_p^*$
is dense and conclude that in this case $\Oh[\frac {1}{l}]^*\subseteq\Oh_{\P}^*$ is dense. For $p=2$ we can choose $l$ such that the closure of $X$ equals
$1+4\Z_2$ and conclude that the closure of $\Oh[\frac{1}{l}]^*$ inside $\Oh_{\P}^*$ equals

\[ \ker(\Oh_{\P}^*\stackrel{\Nr}{\longrightarrow}\Z_2^*\longrightarrow \Z_2^*/(1+4\Z_2)\simeq\{ \pm 1 \} ),\]

c.f. Remark \ref{effective1}. In the special case in which $D$ is the endomorphism algebra of a super-singular
elliptic curve in characteristic $p$, i.e. $\inv_v(D)=0$ for all $v\neq p,\infty$,
this result has been established by different means in \cite[Theorem 0.1]{bl}.\\
2) See Theorem \ref{commapprox} in subsection \ref{cad-1} for a further discussion of the closure of $X\subseteq\Oh_{k,\p}^*$.

\end{example}

The rest of this subsection is devoted to the proof of Theorem \ref{approx}.\\
Remember the groups $G,G'/ \Spec(\Oh_k)$ defined by $G(R)=(\Oh\otimes_{\Oh_k} R)^*$ and $G'(R)=\{ g\in G(R)\, | \,  \Nr(g)=1\}$.

\begin{prop}\label{density} 
The subgroup $G'(\OSfin)\subseteq G'(\Oh_{k,\p})$ is dense.
\end{prop}

\begin{proof}
First note that $G'/\Spec(\Oh_k)$ is representable by an affine
group scheme, hence the injectivity of the homomorphism $G'(\OSfin)\longrightarrow G'(\Oh_{k,\p})$ follows from the injectivity of $\OSfin\hookrightarrow\Oh_{k,\p}$. Secondly, $G'(\Oh_{k,\p})$ is canonically a topological group \cite[Chapter I]{weil} and 
we claim density with respect to this topology. We know that $G'_k:=G'\times_{\Spec(\Oh_k)} \Spec(k)=\Sl_1(D)$ \cite[Section 2.3]{pr} is an inner form of $\Sl_{d,k}$ and thus is $k$-simple, semi-simple
and simply connected. Furthermore, $G_k'\times_{\Spec(k)}\Spec(k_{v_0})=\Sl_n(\tilde{D})$ for
some central skew-field $\tilde{D}$ over $k_{v_0}$ and some $n\ge 1$. Since 
$D_{v_0}$ is not a skew-field by assumption, we have $n\ge 2$ and $\rk_{k_{v_0}}( G'_k\times_{\Spec(k)} \Spec(k_{v_0}))=n-1\ge 1$ \cite[Proposition 2.12]{pr}, i.e. $G'_k$ is isotropic at $v_0$. From strong
approximation \cite[Theorem 5.1.8]{springer} we conclude that

\begin{equation}\label{dense}
G'(k)\cdot G'(k_{v_0})\subseteq G'(\A_k)\mbox{ is dense,}
\end{equation}
where $\A_k$ denotes the ad\`ele-ring of $k$. Fix $x\in G'(\Oh_{k,\p})$
and an open subgroup $U_{\p}\subseteq G'(\Oh_{k,\p})$. Denote by $\tilde{x}\in G'(\A_k)$
the ad\`ele having $\p$-component $x$ and all other components equal to $1$. 
Then

\[ U:=U_{\p}\times\prod\limits_{v\neq\p\mbox{ {\tiny finite}}} G'(\Oh_{k,v})\times\prod\limits_{v\mbox{ {\tiny infinite}}} G'(k_v)\subseteq G'(\A_k) \]

is an open subgroup and by (\ref{dense}) there exist $\gamma\in G'(k)$ and
$\delta\in G'(k_{v_0})$ such that $\gamma\delta\in \tilde{x} U$. Since
$\p\neq v_0$ this implies that $\gamma_{\p}\in\tilde{x}_{\p} U_{\p}= x U_{\p}$, where $\gamma_{\p}$
is the $\p$-component of the principal ad\`ele $\gamma$, equivalently, the image of 
$\gamma$ under the inclusion $G'(k)\subseteq G'(k_{\p})$.
Since $x$ and $U_{\p}$ are arbitrary, we will be done if we can show that
$\gamma\in G'(\OSfin)\subseteq G'(k)$, i.e. that for every finite place $v\not\in\Srm$
we have $\gamma_v\in G'(\Oh_{k,v})$. For $v=\p$ this is clear since $xU_{\p}\subseteq G'(\Oh_{k,\p})$ whereas for $v\neq\p$ we have, using that $\delta_v=1$ since $v\neq v_0\in\Srm$, 
\[ (\gamma\delta)_v=\gamma_v\in(\tilde{x}U)_v=\tilde{x}_v\cdot G'(\Oh_{k,v})= G'(\Oh_{k,v}). \]

\end{proof}

To proceed, we apply (\ref{sequence}) to the inclusion $\OSfin\hookrightarrow\Oh_{k,\p}^*$
to obtain a commutative diagram

\begin{equation}\label{diag}
\xymatrix{  1 \ar[r] & G'(\OSfin) \ar[r] \ar@{^{(}->}[d] & G(\OSfin) \ar@{^{(}->}[d]^{\iota} \ar[r]^{\Nr} & \OSfin^* \ar@{^{(}->}[d] \\
1 \ar[r] & G'(\Oh_{k,\p}) \ar[r] & G(\Oh_{k,\p}) \ar[r]^{\Nr_{\p}} & \Oh_{k,\p}^*. }
\end{equation}

By definition, $G(\OSfin)=(\Oh\otimes_{\Oh_k}\OSfin)^*$ and $G(\Oh_{k,\p})=(\Oh\otimes_{\Oh_k}\Oh_{k,\p})^*=\Oh_{\P}^*$ \cite[Kapitel VI, \S 11, Satz 6]{deu}, so Theorem \ref{approx} is concerned 
with the closure of the image of $\iota$. Recall the subgroup $X\subseteq\OSfin^*$ from (\ref{Xdef}).

\begin{prop}\label{imN} In (\ref{diag}) we have $\im(\Nr)=X\subseteq\OSfin^*$.
\end{prop}

\begin{proof} Eichler's norm theorem \cite[Theorem 1.13]{pr} states that 

\begin{equation}\label{eichler}
\im(\Nr_k:G(k)\longrightarrow k^*)=\{ \alpha\in k^*\, | \, v\in\Sigma_k\mbox{ infinite and }\inv_v(D)\neq 0\mbox{ imply }v(\alpha)>0\},
\end{equation}

and the inclusion $\im(\Nr)\subseteq X$ is trivial by the definition of $X$.\\
From the cohomology sequence associated with (\ref{sequence}) we have

\[ \xymatrix{ G(\OSfin) \ar[r]^{\Nr} \ar@{^{(}->}[d] & \OSfin^* \ar@{^{(}->}[d]
\ar[r] & \HH^1(\Spec(\OSfin),G')\ar[d]^{\iota} \\ 
G(k) \ar[r] & k^*\ar[r] & \HH^1(\Spec(k),G').} \]

Now, $G'/\Spec(\Oh_k)$ is smooth with connected fibers by Theorem \ref{smooth}, and the generic fiber $G'_k$ is an inner form of $\Sl_d$ and is thus $k$-simple, semi-simple and simply connected. Finally, the place $v_0$ is outside $U:=\Spec(\OSfin)$ and since $D_{v_0}$ is not a skew-field, $G'_k$ is isotropic at $v_0$ \cite[Proposition 2.12]{pr}. We can thus apply Proposition \ref{kertrivial} to $G'/U$ to conclude that $\iota$ has trivial kernel. This, jointly with (\ref{eichler}), implies that $X\subseteq\im(\Nr)$.
\end{proof}

We know that $\HH^1(\Spec(\Oh_{k,\p}),G')=0$ from the fact that $G'\times_{\Spec(\Oh_k)}\Spec(\Oh_{k,\p})/\Spec(\Oh_{k,\p})$ is smooth with connected fibers and Lang's Theorem. Hence, in (\ref{diag}), $N_{\p}$ is surjective and we can, using Proposition 
\ref{imN}, rewrite (\ref{diag}) as

\begin{equation}\label{diag2}
\xymatrix{  1 \ar[r] & G'(\OSfin) \ar[r] \ar@{^{(}->}[d]^{\alpha} & (\Oh\otimes_{\Oh_k}\OSfin)^* \ar@{^{(}->}[d]^{\iota} \ar[r]^-{\Nr} & X \ar@{^{(}->}[d]\ar[r] & 1 \\
1 \ar[r] & G'(\Oh_{k,\p}) \ar[r] & \Oh_{\P}^* \ar[r]^{\Nr_{\p}} & \Oh_{k,\p}^* \ar[r] & 1.}
\end{equation}

Since the image of $\alpha$ is dense by Proposition \ref{density} and 
$\Oh_{\P}^*$ is compact, all that remains to be done to conclude the proof
of Theorem \ref{approx} is to apply Proposition \ref{closure} below to (\ref{diag2}).\\
For a subset $Y$ of a topological space $X$ we denote by $\overline{Y}^X$ the 
closure of $Y$ in $X$.

\begin{prop}\label{closure}
Let
\[ \xymatrix{  1 \ar[r] & H' \ar[r] \ar@{^{(}->}[d] & H \ar@{^{(}->}[d] \ar[r]^{\rho} & H'' \ar@{^{(}->}[d]\ar[r] & 1 \\
1 \ar[r] & G' \ar[r] & G \ar[r]^{\pi} & G'' \ar[r] & 1} \] 
be a commutative diagram of first countable topological groups with exact rows, $G$ compact, and such that $H'\subseteq G'$ is dense. Then
\[ \overline{H}^G=\pi^{-1}(\overline{H''}^{G''}).\]
\end{prop}

\begin{proof}
Assume that $g\in\overline{H}^G$. Then $g=\lim\limits_n h_n$ for suitable $h_n\in H$
and $\pi(g)=\lim\limits_n \pi(h_n)\in\overline{H''}^{G''}$. Conversely, given
$g\in G$ with $\pi(g)=\lim\limits_n h_n''$ for suitable $h_n'' \in H''$, choose
$h_n\in H$ with $\rho(h_n)=h_n''$. The sequence $(h_ng^{-1})_n$ in $G$ has a convergent
subsequence, $\tilde{g}:=\lim\limits_i h_{n_i}g^{-1}\in G$. Then $\pi(\tilde{g})=1$, i.e. $\tilde{g}\in G'$
and we have $\tilde{g}=\lim\limits_i h_i'$ for suitable $h_i' \in H'$. The sequence $((h_i')^{-1}
h_{n_i})_i$ in $H$ satisfies $\lim\limits_i (h_i')^{-1} h_{n_i}=\tilde{g}^{-1}\tilde{g}g=g$, hence
$g\in\overline{H}^G$.

\end{proof}

\subsection{Type $^2A_{d-1}$}\label{2ad-1}

Let $D$ be a finite-dimensional skew-field of reduced dimension $d>1$ over $\Q$ carrying a positive involution of the second kind $*$ and assume that $\Oh\subseteq D$ is a maximal order which is stable under $*$. In subsection \ref{geo2ad-1} we associated with these data group schemes $\SU\subseteq\U$ and $\T'\subseteq \T$ over $\Spec(\Oh_k)$ and an open subscheme ${\mathcal U}\subseteq\Spec(\Oh_k)$.\\
To formulate our approximation problem, we fix a prime $0\neq\p\subseteq\Oh_k$ and a finite set of finite places $\Srm$ of $k$ such that
\[ \p\not\in\Srm,\Srm\mbox{ contains all primes of residue characteristic }2\mbox{ ramified in }L/k\mbox{ and}\]
\[ \Srm\mbox{ contains a place }v_0\mbox{ split in }L/k\mbox{ such that for }w_0|v_0\mbox{ }D_{w_0}\mbox{ is not a skew-field.}\]

This implies in particular that $\Spec(\Oh_{k,\Srm})\subseteq {\mathcal U}$. Note that we do not really restrict generality by insisting that $\Srm$ consists of finite place because, unlike the case treated in subsection \ref{ad-1}, the group $\SU$
is anisotropic at every infinite place of $k$.

\begin{theorem}\label{approxinv} 
The closure of $\U(\Oh_{k,\Srm})\subseteq \U(\Oh_{k,\p})$ equals
\[ \{ g\in \U(\Oh_{k,\p})\, | \, \Nr_{\p}(g)\mbox{ lies in the closure of }\T'(\Oh_{k,\Srm})
\subseteq \T'(\Oh_{k,\p})\}.\]
\end{theorem}

See Corollary \ref{smalltorus} for the computation of the closure of $\T'(\Oh_{k,\Srm})\subseteq \T'(\Oh_{k,\p})$ in a special case.\\
Note that \[ \U(\Oh_{k,\Srm})=\{ g\in (\Oh\otimes_{\Oh_k} \Oh_{k,\Srm})^*\, | \, ^*gg=1\} \]
by definition but the structure of $\U(\Oh_{k,\p})$ depends on the splitting behavior of
$\p$ in $L$, c.f. the proof of Theorem \ref{geoofUandSU}. In particular, if 
$\p$ splits in $L/k$ then 
\[ G(\Oh_{k,\p})\simeq\Oh_{D_{\qf}}^*\simeq \Oh_{D_{\qf'}}^* \]
where $\qf$ and $\qf'$ are the primes of $L$ lying above $k$.\\
In the rest of this subsection we give the proof of Theorem \ref{approxinv}.\\

\begin{prop}\label{density2} $\SU(\Oh_{k,\Srm})\subseteq\SU(\Oh_{k,\p})$ is a dense subgroup.
\end{prop}

\begin{proof}
Since $\SU_k:=\SU\times_{\Spec(\Oh_k)}\Spec(k)$
is an outer form of $\Sl_{d,k}$, it is $k$-simple, semi-simple and simply connected. 
In the proof of Theorem \ref{geoofUandSU} we saw that $\SU_{k_{v_0}}\simeq\Sl_1(D_{w_0})$ and since $D_{w_0}$ is not a skew-field by assumption, $\SU_k$
is isotropic at $v_0$.
Now one proceeds as in the proof of 
Proposition \ref{density}.
\end{proof}

Recall from Theorem \ref{geoofUandSU} that we have an exact sequence

\begin{equation}\label{zeile}
1\longrightarrow \SU_{{\mathcal U}}\longrightarrow \U_{{\mathcal U}}\stackrel{\Nr_{{\mathcal U}}}{\longrightarrow} T'_{{\mathcal U}}\longrightarrow1
\end{equation}

over ${\mathcal U}\supseteq\Spec(\Oh_{k,\Srm})$.

\begin{prop}\label{conditions} The diagram obtained by applying (\ref{zeile}) to $\Oh_{k,\Srm}\hookrightarrow\Oh_{k,\p}$ 
\begin{equation}\label{diagram3}
\xymatrix{ 1\ar[r] & \SU(\Oh_{k,\Srm}) \ar[r] \ar@{^{(}->}[d] & \U(\Oh_{k,\Srm}) \ar[r]^{\Nr} \ar@{^{(}->}[d] & \T'(\Oh_{k,\Srm}) \ar[r] \ar@{^{(}->}[d] & 1\\
1\ar[r] & \SU(\Oh_{k,\p}) \ar[r] & \U(\Oh_{k,\p}) \ar[r]^{\Nr_{\p}} & \T'(\Oh_{k,\p}) \ar[r] & 1}
\end{equation}

fulfills the assumptions of Proposition \ref{closure}.
\end{prop}
This finishes the proof of Theorem \ref{approxinv} by applying Proposition \ref{closure} to (\ref{diagram3}).

\begin{proofof} Proposition \ref{conditions}.
Clearly, diagram (\ref{diagram3}) is made up of first-countable groups and is commutative, $\SU(\Oh_{k,\Srm})\subseteq\SU(\Oh_{k,\p})$ is dense by Proposition \ref{density2} and $\U(\Oh_{k,\p})$ is compact. It remains to prove the exactness of the rows, i.e. the surjectivity of $\Nr$ and $\Nr_{\p}$. Since $\SU_{\p}/\Spec(\Oh_{k,\p})$ is smooth with connected fibers by Theorem \ref{geoofUandSU}, Lang's Theorem implies that $\HH^1(\Spec(\Oh_{k,\p}),\SU_{\p})=0$ and thus
$\Nr_{\p}$ is surjective. We now show that $\Nr_k:\U(k)\longrightarrow\T'(k)$ is surjective: we have a commutative diagram with exact rows 

\[ \xymatrix{ \U(k) \ar[r]^{\Nr_k} \ar[d] & \T'(k) \ar[r] \ar[d] & \HH^1(\Spec(k),\SU) \ar[d]^{\simeq} \\
\prod\limits_{v\in\Sigma_k^{\infty}} \U(k_v) \ar[r]^{\prod \Nr_v} & \prod\limits_{v\in\Sigma_k^{\infty}} \T'(k_v) \ar[r] & \prod\limits_{v\in\Sigma_k^{\infty}} \HH^1(\Spec(k_v),\SU). }\]

Here, $\Sigma_k^{\infty}$ denotes the set of infinite places of $k$ and the right-most vertical arrow is an isomorphism by the Hasse-principle for
$\SU\times_{\Spec(\Oh_k)} \Spec(k) $ \cite[Theorem 6.6]{pr}. Hence the surjectivity of $\Nr_k$ will follow from the surjectivity of $\Nr_v$ for all $v\in\Sigma_k^{\infty}$
which is easy to see:\\
For $v\in\Sigma_k^{\infty}$ we have, using \cite[Step IV on page 199]{mumford},

\[ \xymatrix{ \U(k_v)\ar[r]^{\Nr_v}\ar[d]^{\simeq} & \T'(k_v)=\T(k_v)\ar[d]^{\simeq}\\
\{ (x_{i,j})\in\Gl_d(\C)\, | \, (\overline{x_{i,j}})(x_{j,i})=1\} \ar[r]^-{\det} & \{ \alpha\in\C^*\, | \, \alpha\overline{\alpha}=1\}, } \]
where a bar denotes complex conjugation, and the lower horizontal arrow is surjective since it is split by $\alpha\mapsto\diag(\alpha,1,\ldots,1)$. Next we look at the commutative diagram with exact rows

\[ \xymatrix{ \U(\Oh_{k,\Srm})\ar[d]\ar[r]^{\Nr} & \T'(\Oh_{k,\Srm})\ar[r]\ar[d] & \HH^1(\Spec(\Oh_{k,\Srm}),\SU)\ar[d]^{\iota}\\
\U(k)\ar[r]^{\Nr_k} & \T'(k)\ar[r] & \HH^1(\Spec(k),\SU) . } \]

We need to see that $\iota$ has trivial kernel for then the desired surjectivity of $\Nr$ will follow from the already proved surjectivity of $\Nr_k$: since $\Spec(\Oh_{k,\Srm})\subseteq {\mathcal U}$ we know that $\SU/\Spec(\Oh_{k,\Srm})$ is smooth with connected fibers from Theorem \ref{geoofUandSU}, $\SU_k$ is $k$-simple, semi-simple and simply connected and the 
place $v_0$ lies outside $\Spec(\Oh_{k,\Srm})$ and $\SU_k$ is isotropic at $v_0$
as explained in the proof of Proposition \ref{density2}. Hence the kernel of $\iota$ is indeed trivial by Proposition \ref{kertrivial}.
\end{proofof}

\section{The commutative case}\label{commutative}

\subsection{ Type $A_{d-1}$}\label{cad-1}\hfill
In subsection \ref{ad-1} the problem of approximating a local 
unit in a maximal order was reduced to a similar problem involving solely numberfields:\\
Let $k$ be a numberfield, $0\neq\p\subseteq\Oh_k$ a prime dividing the 
rational prime $p$ and $\Sigma$ a possibly empty set of real places of $k$.
For a finite set of finite places $\Srm$ of $k$ not containing $\p$ we consider

\[ X_{\Srm}:=\{x\in\Oh_{k,\Srm}^*\, | \, v(x)>0\mbox{ for all }v\in \Sigma \}\subseteq\Oh_{k,\Srm}^* \]
and wish to understand when $X_{\Srm}\subseteq\Oh_{k,\p}^*=:U_{\p}$ is a dense subgroup. The principal units 
\[ U_{\p}^{(1)}:=1+\p\Oh_{k,\p}\subseteq U_{\p} \]
are canonically a finitely generated $\Z_p$-module and $U_{\p}/U_{\p}^{(1)^p}$ is a finite abelian group.\\
It follows from Nakayama's lemma that a subgroup $Y\subseteq U_{\p}$ 
is dense if and only if the composition $Y\hookrightarrow U_{\p}\longrightarrow U_{\p}/U_{\p}^{(1)^p}$ is surjective: since $U_{\p}$ is pro-finite, $Y\subseteq U_{\p}$ is dense if and only if it surjects onto every finite quotient of $U_{\p}$. Assume that $Y$ does surject onto $U_{\p}/U_{\p}^{(1)^p}$ and $V\subseteq U_{\p}$
is arbitrary of finite index. In order to see that $Y$ surjects onto $U_{\p}/V$ we can assume that $V\subseteq U_{\p}^{(1)^p}$. Then the image of $Y$ in
$U_{\p}/V=\mu_{q-1}\times U_{\p}^{(1)}/V$ surjects onto $\mu_{q-1}$ and $U_{\p}^{(1)}/V$ is a finitely generated $\Z_p$-module which modulo $p$ is generated by the image of $Y$. By Nakayama's lemma, $Y$ surjects onto $U_{\p}/V$.\\
For an 
infinite place $v$ of $k$ we write $k_v^{*,+}$ for the connected component of $1$ inside $k_v^*$, i.e. $k_v^{*,+}\simeq\R^+$ (resp. $k_v^{*,+}\simeq\C^*$) if $v$ is real (resp. complex).\\
We denote by
\[ E^+:=\ker(\Oh_k^*\stackrel{\diag}{\hookrightarrow} \bigoplus\limits_{v\in \Sigma} k_v^*\longrightarrow \bigoplus\limits_{v\in \Sigma} k_v^*/k_v^{*,+} ) \]
the group of global units of $k$ which are positive at all places in $\Sigma$. We write \[ \psi:E^+\subseteq\Oh_k^*\hookrightarrow U_{\p}\] for the inclusion. Then $U_{\p}/\psi(E^+)U_{\p}^{(1)^p}$ is a finite abelian group the minimal number of generators of which
we denote by $g(\p,\Sigma)$.

\begin{theorem}\label{commapprox} In the above situation:\\
i) If $X_{\Srm}\subseteq U_{\p}$ is dense then $|\Srm|\ge g(\p,\Sigma)$.\\
ii) Given a set $\mathrm{T}$ of places of $k$ of density 1, there exists $\Srm$ as above
such that $X_{\Srm}\subseteq U_{\p}$ is dense, $|\Srm|=g(\p,\Sigma)$ and $\Srm\subseteq \mathrm{T}$.\\
iii) \[ g(\p,\Sigma)\leq\left\{ \begin{array}{ccc} [ k_{\p}:\Q_p ] & , & \mbox{ if  }\, \mu_{p^{\infty}}(k_{\p})=\{ 1 \}, \\  1+[k_{\p}:\Q_p] & , & \mbox{ if  }\, \mu_{p^{\infty}}(k_{\p})\neq \{ 1 \}.\end{array} \right. \]
\end{theorem}

\begin{remark} 1) In general, the inequalities in {\sl iii)} are strict: for $k=\Q(\sqrt{2})$, $\p$
dividing $7$ and $\Sigma=\emptyset$ one can check that $g(\p,\Sigma)=0$, i.e. 
$\Oh_k^*\subseteq U_{\p}$ is dense.\\
2) The proof of Theorem \ref{commapprox},{\sl ii)} is rather constructive: one has to find principal prime ideals $(\lambda)$ of $\Oh_k$ with $\lambda$ positive at all places
in $\Sigma$ (this corresponds to being trivial in $\Gal(M/k)$ in the notation of the proof) and determine the image of $\lambda$ in $U_{\p}/\psi(E^+)U_{\p}^{(1)^p}$.
\end{remark}

\begin{proof} We consider the following subgroups of $I_k$, the id\`eles of $k$:
\[ U_K:=\prod\limits_{v\nmid \infty ,v\neq\p} U_v\times U_{\p}^{(1)^p}\times \prod\limits_{v\in \Sigma} k_v^{*,+}\times\prod\limits_{v|\infty, v\not\in \Sigma} k_v^*,\]
\[ U_M:= \prod\limits_{v\nmid \infty} U_v \times \prod\limits_{v\in \Sigma} k_v^{*,+}\times\prod\limits_{v|\infty, v\not\in \Sigma} k_v^*\mbox{   and}\]
\[ U_+:= \prod\limits_{v\nmid \infty} U_v \times \prod\limits_{v| \infty} k_v^{*,+}.\]

Then $U_K\subseteq U_M$ and $k^* U_K\subseteq I_k$ is of finite index. Class field theory, e.g. \cite[Chapter VI]{neukirch}, yields finite abelian extensions $k\subseteq M\subseteq K$ and the upper part of diagram (\ref{cft}) below. The field corresponding to $k^*U_+$ is the big Hilbert class field of $k$ which we denote by $H^+$. Since $k^*U_K\cdot k^*U_+=k^*U_M$ we have $H^+\cap K=M$ and we put $L:=H^+K$.
We have the following diagram of fields

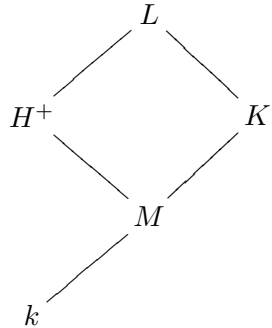
\[ \xymatrix{ & L & \\ H^+ \ar@{-}[ur] & & K \ar@{-}[ul] \\ & M \ar@{-}[ur] \ar@{-}[ul] & \\ k \ar@{-}[ur] & & } \]

and some of the occurring Galois groups are identified as follows:

\begin{equation}\label{cft}
\xymatrix{ 1 \ar[r] & \Gal(K/M) \ar[r]^-{\iota} & \Gal(K/k) \ar[r]^-{\pi} & \Gal(M/k) \ar[r] & 1 \\
1 \ar[r] & k^*U_M/k^*U_K \ar[u]^{\beta}_{\simeq} \ar[r] & I_k/k^*U_K \ar[u]_{\simeq} \ar[r] & I_k/k^*U_M \ar[r] \ar[u]_{\simeq} & 1\\
 & U_{\p}/\psi(E^+)U_{\p}^{(1)^p} \ar[u]^{\alpha}_{\simeq} }
\end{equation}

The isomorphism $\alpha$ is induced by the inclusion $U_{\p}\hookrightarrow k^*U_M$: one has $k^*U_M=k^*U_{\p}U_K$, hence 
\[ \xymatrix{ k^*U_M/k^*U_K=k^*U_{\p}U_K/k^*U_K & U_{\p}/(U_{\p}\cap k^*U_K) \ar[l]_-{\simeq} } \]
and $U_{\p}\cap k^*U_K=k^*U_{\p}\cap U_K=\psi(E^+)U_{\p}^{(1)^p}$.\\
To prove {\sl i)}, assume that $X_{\Srm}\subseteq U_{\p}$ is dense. Then $X_{\Srm}\subseteq U_{\p}\longrightarrow U_{\p}/U_{\p}^{(1)^p}$ is surjective, hence so is $X_{\Srm}/E^+\longrightarrow U_{\p}/\psi(E^+)U_{\p}^{(1)^p}$. 
The group $X_{\Srm}/E^+$ is easily seen to be torsion-free and Dirichlet's
unit Theorem determines its rank, hence $X_{\Srm}/E^+\simeq\Z^{|\Srm|}$ and 
$|\Srm|\ge g(\p,\Sigma)$.\\
To prove {\sl ii)}, fix generators $x_i\in U_{\p}/\psi(E^+)U_{\p}^{(1)^p}$ ($1\leq i\leq g(\p,\Sigma)$). Let $\sigma_i\in\Gal(L/M)\subseteq\Gal(L/k)$ be the unique element such that
$\sigma_i|_{H^+}=\id$ and $\sigma_i|_K=(\iota\beta\alpha)(x_i)$. Note that $(\iota\beta\alpha)(x_i)|_M=(\pi\iota\beta\alpha)(x_i)=\id$ by (\ref{cft}). By Chebotarev's
density Theorem \cite[Chapter VII, Theorem 13.4]{neukirch}, there is a finite place 
$v_i\in \mathrm{T}$, unramified in $L/k$ such that $\sigma_i=\Frob_{v_i}^{-1}$, where $\Frob_{v_i}$ denotes the Frobenius at the place $v_i$, in $\Gal(L/k)$. Then $(\iota\beta\alpha)(x_i)=\Frob_{v_i}^{-1}$ in $\Gal(K/k)$.
Since $\Frob_{v_i}|_{H^+}=\sigma_i^{-1}|_{H^+}=\id$, the prime ideal $\p_i\subseteq
\Oh_k$ corresponding to $v_i$ is principal, generated by a totally positive
element $\pi_i\in\Oh_k$ \cite[Chapter VI, Theorem 7.3]{neukirch}.
We claim that the image of $\pi_i$ in $U_{\p}/\psi(E^+)U_{\p}^{(1)^p}$ equals $x_i$:\\
To see this, we apply the Artin-map $(-,K/k):I_k\longrightarrow\Gal(K/k)$ to the
identity $\pi_i=\pi_{i,\p}\cdot(\frac{\pi_i}{\pi_{i,\p}})$ in $I_k$, where $\pi_{i,\p}$
denotes the id\`ele having $\pi_i$ as its $\p$-component and all other components equal to $1$. By Artin-reciprocity we obtain $1=(\pi_{i,\p},K/k)(\frac{\pi_i}{\pi_{i,\p}},K/k)$. Denoting $y:= \frac{\pi_i}{\pi_{i,\p}}$ we have $(y,K/k)=\prod\limits_v (y_v,K_v/k_v)$ \cite[Chapter VI, Theorem 5.6]{neukirch} and evaluate the local terms $(y_v,K_v/k_v)$
as follows:\\
For $v=\p$ we obtain $1$ since $y_{\p}=1$; for $v\neq \p,v_i$ finite we obtain
$1$ since $y_v\in\Oh_{k,v}^*$ and $v$ is unramified in $K/k$; for $v=v_i$ we obtain
$\Frob_{v_i}$ since $K/k$ is unramified at $v_i$ and $y_{v_i}\in\Oh_{k,v_i}$
is a local uniformizer; finally, for $v|\infty$ we obtain $1$ since $y_v>0$ because $\pi_i$
is totally positive.\\
Hence $(\pi_{i,\p},K/k)=\Frob_{v_i}^{-1}=(\iota\beta\alpha)(x_i)$ in $\Gal(K/k)$.
Denoting by $\tau: U_{\p}\longrightarrow U_{\p}/\psi(E^+)U_{\p}^{(1)^p}$ the projection we have
$(\pi_{i,\p},K/k)=(\iota\beta\alpha\tau)(\pi_{i,\p})$ by construction, hence 
$x_i=\tau(\pi_{i,\p})$ by the injectivity of $\iota\beta\alpha$. This establishes the
above claim saying that the global elements $\pi_i\in\Oh_k$ have the prescribed image $x_i$ in $U_{\p}/\psi(E^+)U_{\p}^{(1)^p}$. To conclude the proof of {\sl ii)},
put $\Srm:=\{ v_i\, |\, 1\leq i\leq g(\p,\Sigma)\} $ and note that $\pi_i\in X_{\Srm}$ with this choice of S, hence
$X_{\Srm}\longrightarrow U_{\p}/\psi(E^+)U_{\p}^{(1)^p}$ is surjective and since $E^+\subseteq X_{\Srm}$, so is $X_{\Srm}\longrightarrow U_{\p}/U_{\p}^{(1)^p}$, i.e. $X_{\Srm}\subseteq U_{\p}$ 
is dense and by construction we have $\Srm\subseteq\mathrm{T}$ and $|\Srm|=g(\p,\Sigma)$.\\
To see {\sl iii)} we use 
\[ U_{\p}=\mu_{q-1}\times U_{\p}^{(1)}\simeq\mu_{q-1}\times\mu_{p^{\infty}}(k_{\p})\times\Z_p^{ [k_{\p}:\Q_p ]} ,\]

where $q=|\Oh_{k,\p}/\p \Oh_{k,\p}|$ \cite[Chapter II, Theorem 5.7, i)]{neukirch} which 
implies that the upper bound claimed in {\sl iii)} is in fact the minimal number of generators of $U_{\p}/U_{\p}^{(1)^p}$ which obviously is greater than or equal to the minimal number of generators of $U_{\p}/\psi(E^+)U_{\p}^{(1)^p}$, i.e. $g(\p,\Sigma)$.
\end{proof}

\subsection{Type $^2A_{d-1}$}\label{c2ad-1}

In Theorem \ref{approxinv} we reduced the problem of approximating a local unit 
of a maximal order carrying a positive involution of the second kind by global {\em unitary} units to an approximation problem
for a specific integral model $\T'$ of a one-dimensional anisotropic torus over a totally real number-field. This approximation problem seems to be substantially harder than the
problem settled in Theorem \ref{commapprox} and we only treat the following 
special case here.\\
Let $k$ be an imaginary quadratic field in which the rational prime $p$ 
splits, $p\Oh_k=\p\overline{\p}$, and put

\[ \T:=\ker(\Res^{\Oh_k}_{\Z}(\G_{m,\Oh_k})\stackrel{\Nr^k_{\Q}}{\longrightarrow}\G_{m,\Z}).\]

\begin{theorem}\label{approxtorus} In the above situation, there exist infinitely many rational 
primes $l\ne p$ which split in $k/\Q$ and are such that $\T(\Z[1/l])\subseteq \T(\Z_p)$ is a dense subgroup.
\end{theorem}

\begin{proof} Note that for every rational prime $l\neq p$

\begin{equation}\label{dingsda}
\T(\Z[1/l])=\{ \alpha\in \Oh_k[1/l]^*\, | \,
\alpha\overline{\alpha}=1\} \subseteq \T(\Z_p)=U_{\p}\simeq\Z_p^*, 
\end{equation}

the local units of $k$ at $\p$, the final equalities following from the fact that $p$ splits in $k$. Here, $-$ denotes complex conjugation. The following proof is similar to the argument of Theorem \ref{commapprox},{\sl ii)}
but extra care is needed to deal with the norm condition
$\alpha\overline{\alpha}=1$.\\
Consider the following subgroups of the id\`eles of $k$:

\[ U_K:=\prod\limits_{v\neq\p, \overline{\p}\,\,\mbox{\tiny finite }} U_v\times 
U_{\p}^{(1)^p} \times U_{\overline{\p}}^{(1)^p}\times \prod\limits_{v |
  \infty}k_v^* \mbox{ and}\]

\[ U_H:= \prod\limits_{v\,\,\mbox{\tiny finite}}U_v\times\prod\limits_{v |
  \infty} k_v^*. \]

We have a corresponding tower of abelian extensions $k\subseteq
H\subseteq K$ and since $U_K$ is stable under $\Gal(k/\Q)$, the
extension $K/\Q$ is Galois, though rarely abelian. We have an
isomorphism

\[ \phi:
U_{\p}U_{\overline{\p}}/U_{\p}^{(1)^p}U_{\overline{\p}}^{(1)^p}\Oh_k^*\simeq\frac{U_{\p}/U_{\p}^{(1)^p}\times
  U_{\overline{\p}}/U_{\overline{\p}}^{(1)^p}}{\Oh_k^*}\stackrel{\simeq}{\longrightarrow}\Gal(K/H)
  \]

induced by the Artin-map, where $\Oh_k^*$ is embedded
diagonally. Since $p$ splits in $k$ we have $U_{\p}\simeq\Z_p^*$, $U_{\p}/U_{\p}^{(1)^p}\Oh_k^*$ is
cyclic and we fix a generator $x$ of this group. By Chebotarev's Theorem applied
to $K/\Q$ there exist infinitely many rational primes $l\neq p$, unramified in
$K/\Q$ and such that for a suitable prime $\Lambda$ of $K$ lying above $l$ we have

\[ \Frob_{\Lambda | l} ^{-1} = \phi([ (x,1) ] )\mbox{ in
  }\Gal(K/H)\subseteq\Gal(K/\Q). \]

We claim that every such $l$ satisfies the conclusion of Theorem
\ref{approxtorus}.\\
Put $\lambda := \Lambda | _k$. Since $ (\Frob_{\Lambda | l})|_ H=\id $,
$l$ is split in $k/\Q$ and
$\lambda$ is a principal ideal of $\Oh_k$ a generator of which we
denote by $\pi$. Then 

\[ \beta := \frac{\pi}{\overline{\pi}} \in \left\{ \alpha\in\Oh_k[1/l]^* \,
    | \, \alpha\overline{\alpha}= 1 \right\}=T(\Z[1/l]) , \]

and we claim that $\beta$ goes to $x$ under the map induced by (\ref{dingsda}): as in the proof of Theorem
\ref{commapprox},{\sl ii)} one sees that

\[ (\pi_{\p},\pi_{\overline{\p}})=[(x,1)]\mbox{ and similarly }\]
\[ ( (\overline{\pi})_{\p}, (\overline{\pi})_{\overline{\p}})=[
(1,x)]\mbox{ in } \frac{U_{\p}/U_{\p}^{(1)^p}\times
  U_{\overline{\p}}/U_{\overline{\p}}^{(1)^p}}{\Oh_k^*}
  , \]

hence indeed \[ (\beta_{\p},\beta_{\overline{\p}})=[ (x,x^{-1}) ] \]
and a fortiori $\beta_{\p}=x$ in $U_{\p}/U_{\p}^{(1)^p}\Oh_k^*$.
Since we have $\Oh_k^*\subseteq \T(\Z[1/l])$ because $\Oh_k^*$ consists of
roots of unity which have norm 1, we are done.

\end{proof}

To use Theorem \ref{approxinv} we must study approximation for the open subgroup-scheme $\T'\subseteq\T$ obtained from $\T$ by removing the non-identity components of finitely many special fibers of $\T$, c.f. subsection \ref{geo2ad-1}. Let $\mu\subseteq\Oh_k^*$ denote the group of roots of unity. While we have $\mu\subseteq\T(\Z)$, and this was used at the end of the proof of Theorem \ref{approxtorus}, in general we also have $-1\not\in\T'(\Z[1/l])$.

\begin{cor}\label{smalltorus} In the above situation there exist
infinitely many rational primes $l\neq p$ which split in $k/\Q$ and are such that the closure of $\T'(\Z[\frac{1}{l}])\subseteq\T'(\Z_p)=\T(\Z_p)$ has index at most $|\mu|$.
\end{cor}

\begin{proof}
We have $\T'\times_{\Spec\Z} \Spec(\Z_p)\stackrel{\simeq}{\longrightarrow}\T\times_{\Spec\Z} \Spec(\Z_p)$ by the construction of $\T'$ and the fact that $p$ if unramified (in fact, split) in $k/\Q$. Now observe that the element $\beta\in\T(\Z[1/l])$ constructed in the
proof of Theorem \ref{approxtorus} satisfies $\beta\in\T'(\Z[1/l])$.
\end{proof}

\begin{cor}\label{unitaryapplication}
Let $D$ be a finite-dimensional skew-field over $\Q$ of reduced dimension $d>1$ with a positive involution
 of the second kind $*$ and $\Oh\subseteq D$ a maximal order, stable under $*$. Assume that the center of $D$ is an imaginary quadratic field $k$ and let $p\neq 2$
be a rational prime which splits in $k$ and $\P\subseteq\Oh$ a prime lying above
$p$. Then there exists a rational prime $l\neq p$ such that the closure of 

\[ \left\{ g\in\Oh\left[ \frac{1}{2l} \right] \mid\,  ^*gg=1\right\}\subseteq\Oh_{\P}^*\]

has index at most $|\mu|$.
\end{cor}

\begin{proof}
From the data $(D,*)$ and $\Oh\subseteq D$ we construct group-schemes $\SU\subseteq\U$ and $\T'\subseteq\T$ over $\Spec(\Z)$ as in subsection \ref{geo2ad-1}. Using Corollary \ref{smalltorus} we choose a prime $l\neq2,p$ which splits in $k/\Q$
such that the closure of $\T'(\Z[1/l])\subseteq\T'(\Z_p)$ has index at most $|\mu|$ and such that for every place $\lambda$ of $k$ lying above $l$ we have
$\inv_{\lambda}(D)=0$. We apply Theorem \ref{approxinv} with $\Srm:=\{ 2,l\}$
to conclude that the index of the closure of $\U(\Z[1/2l])\subseteq\U(\Z_p)$
equals the index of the closure of $\T'(\Z[1/2l])\subseteq\T'(\Z_p)$ and is thus bounded above by $|\mu|$. It remains to recall that

\[ \U\left(\Z\left[ \frac{1}{2l}\right]\right)=\left\{ g\in\Oh\left[\frac{1}{2l}\right]\, | \, ^*gg=1 \right\} \]
and, since $p$ splits in $k/\Q$,

\[ \U(\Z_p)\simeq\Oh_{\P}^*.\]
\end{proof}

\begin{remark} The conclusion of Corollary \ref{unitaryapplication} can be sharpened in special cases, for example: if the reduced dimension of $D$ is odd and $2$ is unramified in $k/\Q$, then ($p=2$ being allowed) there is a rational prime $l\neq p$ such that $\left\{ \alpha\in\Oh[\frac{1}{l}]\mid\,  ^*\alpha\alpha=1\right\}\subseteq\Oh_{\P}^*$ is dense. This is because $d$ being odd implies that $\T'=\T$ and $2$ being unramified implies that ${\mathcal U}=\Spec(\Z)$.
\end{remark}

\section{Applications} \label{applications}

\subsection{Extending automorphisms of $p$-divisible
  groups} \label{extendaut}

Here we explain the application of some of the above results to the following problem:\\
Let $k$ be a finite field of characteristic $p$ and $A/k$ a simple
abelian variety such that $\End_k(A)$ is a maximal order in the
skew-field
$D:=\End_k(A)\otimes_{\Z}\Q$. The center $K$ of $D$ is a numberfield
and $K\cap \End_k(A)=\Oh_K$ is its ring of integers.\\
The $p$-divisible group of $A/k$ \cite{tatepdiv} splits as

\begin{equation}\label{split1}
A[p^{\infty}]=\prod\limits_{\p | p} A[\p^{\infty}],
\end{equation}

the product extending over all primes $\p$ of $\Oh_K$ dividing $p$. According
to J. Tate, c.f. \cite[Theorem 6]{milnewaterhouse}, the canonical homomorphism

\begin{equation}\label{tateiso}
 \End_k(A)\otimes_{\Z}\Z_p\stackrel{\simeq}{\longrightarrow}\End_k(A[p^{\infty}])
\end{equation}

is an isomorphism. We have 

\[ \End_k(A)\otimes_{\Z}\Z_p\simeq \prod\limits_{\p | p}
\End_k(A)\otimes_{\Oh_K}\Oh_{K,\p}\simeq \prod\limits_{\p | p}
\End_k(A)_{\P} \]

with $\P$ the unique prime of $\End_k(A)$ lying above $\p$. Similarly,
(\ref{split1})
implies that 

\[ \End_k(A[p^{\infty}])\simeq \prod\limits_{\p | p}
\End_k(A[\p^{\infty}]). \]

These decompositions are compatible with (\ref{tateiso}) in that the
canonical homomorphism 

\[ \End_k(A)\otimes_{\Oh_K}\Oh_{K,\p}\stackrel{\simeq}{\longrightarrow}
\End_k(A[\p^{\infty}]) \]

is an isomorphism for every $\p | p$. We fix some $\p | p$ and ask for
a finite set S of finite primes of $K$ such that $\p\notin\Srm$ and

\begin{equation}\label{inclusion}
\left(\End_k(A)\otimes_{\Oh_K}\Oh_{K,S}\right)^*\hookrightarrow\Aut_k(A[\p^{\infty}])
\end{equation}

is a dense subgroup. Note that this density is equivalent to the
following assertion:\\
For every $\alpha\in\Aut_k(A[\p^{\infty}])$ and integer $n\ge 1$ there is an
isogeny $\phi\in\End_k(A)$ of degree divisible by primes in S only and
some
$x\in\Oh_{K,S}^*$ such that \[ \phi x|_{A[\p^{n}]}=\alpha
|_{A[\p^{n}]} ,\]

i.e. the quasi-isogeny $\phi x$ of $A$ extends the truncation at
arbitrary finite level $n$ of $\alpha$.\\
By Theorem \ref{approx}, the inclusion (\ref{inclusion}) is dense if
and only if
$X\subseteq U_{\p}$ is dense where $X\subseteq\Oh_K^*$ is the subgroup
of global units which are positive at all real places of $K$ at which $D$
does not split and $U_{\p}:=\Oh_{K,\p}^*$ are the local units of $K$ at $\p$. 
The density of $X\subseteq U_{\p}$ in turn is firmly controlled by
Theorem \ref{commapprox}. We would like to illustrate all of this with
some examples:\\
According to the Albert-classification \cite[Theorem 2, p. 201]{mumford},
note
that types I and II do not occur over finite fields, there
are two possibilities:\\
{\bf Type III}: here, $K$ is a totally real numberfield and
$D/K$
is a totally definite quaternion algebra. The simplest such case
occurs if $A/k$ is a super-singular elliptic curve with
$\End_k(A)=\End_{\overline{k}}(A)$. In this case, it follows from
Example \ref{tralala},2) that, in case the characteristic of $k$ is different from $2$, for a suitable prime $l$ 

\[ \left( \End_k(A)\left[\frac{1}{l}\right] \right)^*\hookrightarrow
\Aut_k(A[p^{\infty}]) \] is dense.\\
To see another example of this type, let $A/\F_p$ correspond to a
$p$-Weil number $\pi$ with $\pi^2=p$. Then $\dim(A)=2$ and
$A\otimes_{\F_p}\F_{p^2}$ is isogeneous to the square of a
super-singular elliptic curve.
We have $K=\Q(\sqrt{p})$ and $\p=(\sqrt{p})\Oh_K$, hence $A[\p^{\infty}]=
A[p^{\infty}]$. Furthermore, $\Oh_K^*=\{ \pm 1 \}\times\epsilon^{\Z}$
for a fundamental unit $\epsilon$ and $X\subseteq\Oh_K^*$ is of index
4.
To find a small set S such that (\ref{inclusion}) is dense one first
needs to compute the minimal number of generators of
$ U_{\p}/XU_{\p}^{(1)^p} $,
denoted $g(\p,\Sigma)$ in Theorem \ref{commapprox} where, in the
present situation, $\Sigma$ consists of both the infinite
places of $K$. For $p=2$ one can choose $\epsilon=1+\sqrt{2}$, then
$X=\epsilon^{2\Z}$. Since $U_{\p}/U_{\p}^{(1)^2}\simeq \F_2^3$ and
$\epsilon^2\notin U_{\p}^{(1)^2}$ one gets $g(\p,\Sigma)=2$.\\
For $p=3$ we may take $\epsilon=2+\sqrt{3}$, then $X=\epsilon^{2\Z}$
again. Since now $U_{\p}/U_{\p}^{(1)^3} = \mu_2\times\F_3^2\simeq
\Z/6\times \Z/3$ the fact that $\epsilon^2\notin U_{\p}^{(1)^3}$ is not
enough to conclude that $g(\p,\Sigma)=1$. However, one checks in
addition\
that $\epsilon^2\in U_{\p}^{(1)}$, and concludes that
$U_{\p}/XU_{\p}^{(1)^3}\simeq\Z/6$ and hence indeed
$g(\p,\Sigma)=1$.\\
For $p\ge 5$ one has $U_{\p}/U_{\p}^{(1)^p}=\mu_{p-1}\times \F_p^2$ and 
since $\mu_{p-1}\not\subseteq K$ the image of a generator of $X$ in 
$U_{\p}/U_{\p}^{(1)^p}$ will have non-trivial projection to $\F_p^2$
and one concludes that $g(\p,\Sigma)=1$.\\

{\bf Type IV:} In this case, $K$ is a CM-field and $X=\Oh_K^*$. The
easiest such example occurs for an ordinary elliptic curve and we give
two examples:\\
A solution of $\pi^2+5=0$ is a $5$-Weil number to which there
corresponds an elliptic curve $E/\F_5$ with $K=D=\Q(\sqrt{5})$. For 
$\p=(\sqrt{5})\Oh_K$ one has $U_{\p}/U_{\p}^{(1)^p}=\mu_4\times \F_5^2$
and since $\Oh_K^*=\{ \pm 1 \}$ one gets $
U_{\p}/XU_{\p}^{(1)^p}\simeq\Z/10\times \Z/5$, hence
$g(\p,\Sigma)=2$.\\
Similarly, a solution of $\pi^2-4\pi+5=0$ gives an elliptic curve over
$\F_5$ with $D=K=\Q(i)$ and since $5$ splits in $K$ one has
$U_{\p}/XU_{\p}^{(1)^p}\simeq\Z/10$, hence $g(\p,\Sigma)=1$ in
this case.\\

Finally, we leave it as an easy exercise to an interested reader to
check that for every prime $p$ and integer $N\ge 1$ there exists a 
simple abelian variety $A/\F_p$ such that every set S for which
(\ref{inclusion}) is dense necessarily satisfies $|S|\ge N$.\\

\subsection{ A dense subgroup of quasi-isogenies in the Morava stabilizer group} \label{qisogdense}

Let $p$ be a prime and $n\ge 1$ an integer. The $n$-th Morava-stabilizer
group $\Smo_n$ is the group of units of the maximal order of the central skew-field over $\Q_p$ of Hasse-invariant $\frac{1}{n}$.\\
In this section we will construct an abelian variety $A/k$ over a finite field
$k$ of characteristic $p$ such that for a suitable prime $l$ the group $(\End_k(A)[\frac{1}{l}])^*$ is canonically a dense subgroup of $\Smo_n$. We will completely ignore the case $n=1$ as it is very well understood. In case $n=2$ one can take for $A$
a super-singular elliptic curve \cite{bl} and the resulting 
dense subgroup of $\Smo_2$ has been used to great advantage in the construction of a modular resolution of the $K(2)$-local sphere \cite{behrens}.\\
For general $n$ we remark that, since $\End_k(A)\otimes_{\Z}\Z_p\simeq\End_k(A[p^{\infty}])$, in order that $\End_k(A)$ have a relation with $\Smo_n$ one needs $A[p^{\infty}]\otimes_k\overline{k}$
to have an isogeny factor of type $G_{1,n-1}$ \cite[IV,\S 2,2.]{manin}. By the symmetry
of $p$-divisible groups of abelian varieties \cite[IV, \S 3, Theorem 4.1]{manin}, there must then
also be a factor of type $G_{n-1,1}$ which shows that $n=2$ is somewhat special
since $(1,n-1)=(n-1,1)$ in this case. For $n\ge 3$ the above considerations
imply that the sought for abelian variety must be of dimension at least $n$, 
as already remarked by D. Ravenel \cite[Corollary 2.4 (ii)]{ravenel}. Following suggestions
of M. Behrens and T. Lawson we will be able to construct $A$ having this minimal possible dimension. We start by constructing a suitable isogeny-class as follows.

\begin{prop} \label{isogclass}
Let $p$ be a prime and $n\ge 3$ an integer. Then there is a simple abelian variety $A/\F_{p^n}$ such that the center of $\, \End_{\F_{p^n}}(A)\otimes_{\Z}\Q\, $ is
an imaginary quadratic field in which $p$ splits into, say, $\p$ and
$\p'$ such that $\inv_{\p}(\End_{\F_{p^n}}(A)\otimes_{\Z}\Q)=\frac{1}{n}$, 
$\inv_{\p'}(\End_{\F_{p^n}}(A)\otimes_{\Z}\Q)=\frac{-1}{n}$ and $\dim(A)=n$.\\
Furthermore, $A$ is geometrically simple with $\End_{\overline{\F}_{p^n}}(A)\otimes_{\Z}\Q=\End_{\F_{p^n}}(A)\otimes_{\Z}\Q.$
\end{prop}

\begin{proof} We use Honda-Tate theory, see \cite{milnewaterhouse} for an exposition. Let $\pi\in\overline{\Q}$ be a root of $f:=x^2-px+p^n\in\Z[x]$. Since the 
discriminant of $f$ is negative, $\pi$ is a $p^n$-Weil number and we choose $A/\F_{p^n}$ simple associated with the conjugacy class of $\pi$. Then $\Q(\pi)$ is
an imaginary quadratic field and is the center of $\End_{\F_{p^n}}(A)\otimes_{\Z}\Q$. Since $n\ge 3$ the Newton polygon of $f$ at $p$ has different slopes $1$ and $n-1$
which shows that $f$ is reducible over $\Q_p$ \cite[Chapter II, Theorem 6.4]{neukirch}, hence $p$ splits in $\Q(\pi)$
into $\p$ and $\p'$ and, exchanging $\pi$ and $\overline{\pi}$ if necessary, 
we can assume that $v_{\p}(\pi)=1$ and $v_{\p}(\overline{\pi})=n-1$. Then \cite[Theorem 8, 4.]{milnewaterhouse}

\[ \inv_{\p}(\End_{\F_{p^n}}(A)\otimes_{\Z}\Q)=\frac{v_{\p}(\pi)}{v_{\p}(p^n)}[\Q(\pi)_{\p}:\Q_p]=\frac{1}{n}\mbox{ and similarly}\]

\[ \inv_{\p'}(\End_{\F_{p^n}}(A)\otimes_{\Z}\Q)=\frac{n-1}{n}=\frac{-1}{n}. \]

Furthermore \cite[Theorem 8, 3.]{milnewaterhouse}, $2\cdot\dim(A)=[\End_{\F_{p^n}}(A)\otimes_{\Z}\Q:\Q(\pi)]^{1/2}\cdot [\Q(\pi):\Q]=2n$. The final statement follows easily from the fact that $\pi^k\not\in\Q$ for all $k\ge 1$, c.f. \cite[Proposition 3(2)]{hz}, which in turn is evident since $v_{\p}(\pi)\neq v_{\p}(\overline{\pi})$.
\end{proof}

Since the properties of $A/\F_{p^n}$ in Proposition \ref{isogclass} are invariant
under $\F_{p^n}$-isogenies, we can, and do, choose $A/\F_{p^n}$ having these properties such that in addition $\End_{\F_{p^n}}(A)\subseteq \End_{\F_{p^n}}(A)\otimes_{\Z}\Q$ is a maximal order \cite[proof of Theorem 3.13]{waterhouse}. Denoting by $\P\subseteq \End_{\F_{p^n}}(A)$ the unique prime lying above the prime $\p$ constructed in Proposition \ref{isogclass}, we have $(\End_{\F_{p^n}}(A))_{\P}^*=\Smo_n$ since $\inv_{\p}(\End_{\F^{p^n}}^0(A)\otimes_{\Z}\Q)=1/n$. We choose a prime $l$ as follows:
If $p\neq 2$ we take $l$ to be a topological generator of $\Z_p^*$. For $p=2$ we take $l=5$. 

\begin{remark}\label{effective1}
Note that for $p\neq 2$ a prime $l\neq p$ topologically generates $\Z_p^*$ if and only if $(l\mbox{ mod }p^2)$ generates $(\Z/p^2)^*$. Hence, by Dirichlet's Theorem on primes in arithmetic
progressions, the set of all such $l$ has a density equal to $((p-1)\phi(p-1))^{-1}>0$ and is thus infinite. Such an $l$ can be found rather effectively: given $l\neq p$, compute $\alpha_k:=(l^{p(p-1)/k}\mbox{ mod }p^2)$ for all primes $k$ dividing $p(p-1)$. If for all $k$, $\alpha_k\not\equiv 1\,(p^2)$, then $l$ is suitable.
\end{remark} 

\begin{theorem}\label{modular1} In the above situation 

\[ (\End_{\F_{p^n}}(A)[ \frac{1}{l}])^*\hookrightarrow (\End_{\F_{p^n}}(A))_{\P}^*=\Smo_n \]
is a dense subgroup.
\end{theorem}

\begin{proof} We apply Theorem \ref{approx} with $\Oh:=\End_{\F_{p^n}}(A)$,
$k:=\Q(\pi)$, $\p$ the prime of $\Oh_k$ constructed in Proposition
\ref{isogclass} and $\Srm:=\{ \infty,l\}$ the set consisting of the unique
infinite place $\infty$ of $k$ and all places dividing $l$. Clearly, $\p\not\in
\Srm$ and $D:=\Oh\otimes_{\Z}\Q$ is not a skew-field at $\infty$ since $k_\infty\simeq\C$ and $n>1$.
Using the notation of Theorem \ref{approx} we have $\OSfin=\Oh_k[1/l]$ and
$X=(\Oh_k[1/l])^*$ since $k$ has no real place. Theorem \ref{approx} shows that the claim of Theorem \ref{modular1} is equivalent to the
density of $(\Oh_k[1/l])^*\subseteq\Oh_{k,\p}^*\simeq\Z_p^*$. Since $l\in (\Oh_k[1/l])^*$, this density is clear for $p\neq 2$ by our choice of $l$ whereas for $p=2$ we have that $\{\pm 1\} \times 5^{\Z}\subseteq\Z_2^*$ is dense and $-1,5\in (\Oh_k[1/5])^*$.
\end{proof}


\begin{thebibliography}{9999}
  
\bibitem[B]{behrens} M. Behrens, A modular description of the $K(2)$-local sphere at the prime 3,  Topology  {\bf 45}  (2006),  no. 2, 343--402.
\bibitem[BL1]{topautforms} M. Behrens, T. Lawson, Topological automorphic forms, http://front.math.ucdavis.edu/math.AT/0702719
\bibitem[BL2]{bl} M. Behrens, T. Lawson, Isogenies of elliptic curves and the Morava stabilizer group, J. Pure Appl. Algebra  {\bf 207}  (2006),  no. 1, 37--49. 
\bibitem[BLR]{blr} S. Bosch, W. L\"utkebohmert, M. Raynaud, N\'eron models, Ergebnisse der Mathematik und ihrer Grenzgebiete (3) {\bf 21}, Springer-Verlag, Berlin, 1990.
\bibitem[BT1]{bt1} F. Bruhat, J. Tits, Groupes r\'eductifs sur un corps local II, Sch\'emas en groupes, Existence d'une donn\'ee radicielle valu\'ee, Inst. Hautes \'Etudes Sci. Publ. Math. No. {\bf 60} (1984), 197--376.
\bibitem[BT2]{bt2} F. Bruhat, J. Tits, Sch\'emas en groupes et immeubles des groupes classiques sur un corps local, Bull. Soc. Math. France  {\bf 112}  (1984),  no. 2, 259--301. 
\bibitem[CTS]{sansuc} J.-L. Colliot-Th\'el\`ene, J.-J.  Sansuc, Principal homogeneous spaces under flasque tori: applications, J. Algebra {\bf 106} (1987), no. 1, 148--205. 
\bibitem[D]{deu} M. Deuring, Algebren, zweite Auflage, Ergebnisse der Mathematik und ihrer Grenzgebiete, Band {\bf 41}, Springer-Verlag, Berlin-New York, 1968.
\bibitem[Gi]{gille} P. Gille, Torseurs sur la droite affine, Transform. Groups {\bf 7} (2002), no. 3, 231--245. 
\bibitem[EGA IV$_3$]{ega43} A. Grothendieck, \'El\'ements de g\'eom\'etrie alg\'ebrique, IV, \'Etude locale des sch\'emas et des morphismes de sch\'emas III, Inst. Hautes \'Etudes Sci. Publ. Math. No. {\bf 28}, 1966.
\bibitem[EGA IV$_4$]{ega44} A. Grothendieck, \'El\'ements de g\'eom\'etrie alg\'ebrique IV, \'Etude locale des sch\'emas et des morphismes de sch\'emas IV, Inst. Hautes \'Etudes Sci. Publ. Math. No. {\bf 32}, 1967.
\bibitem[GHMR]{hennetc} P. Goerss, H.-W. Henn, M. Mahowald, C. Rezk, A resolution of the $K(2)$-local sphere at the prime 3, Ann. of Math. (2) {\bf 162} (2005), no. 2, 777--822.
\bibitem[H]{henn2} H.-W. Henn, On finite resolutions of $K(n)$-local spheres, available at: http://hopf.math.purdue.edu/cgi-bin/generate?/Henn/kn-res-ded
\bibitem[Hi]{hida} H. Hida, $p$-adic automorphic forms on Shimura varieties, Springer Monographs in Mathematics, Springer-Verlag, New York, 2004.
\bibitem[HS]{hoveystrickland} M. Hovey, N. Strickland, Morava $K$-theories and localisation, Mem. Amer. Math. Soc. {\bf 139} (1999), no. 666.
\bibitem[HZ]{hz} E. Howe, H. Zhu, On the existence of absolutely simple abelian varieties of a given dimension over an arbitrary field, J. Number Theory {\bf 92} (2002), no. 1, 139--163.
\bibitem[K]{kleinert} E. Kleinert, Units in skew fields, Progress in Mathematics, {\bf 186}, Birkh\"auser Verlag, Basel, 2000.
\bibitem[KMRT]{bookinv} M. Knus, A. Merkurjev, M. Rost, J-P. Tignol, The Book of Involutions, American Mathematical Society Colloquium Publications, {\bf 44}, American Mathematical Society, Providence, RI, 1998.
\bibitem[Ma]{manin} J. Manin, Theory of commutative formal groups over fields of finite characteristic, Uspehi Mat. Nauk {\bf 18} (1963), no. 6 (114), 3--90.
\bibitem[Mi]{milneetlae} J. Milne, \'Etale cohomology, Princeton Mathematical Series, {\bf 33}, Princeton University Press, Princeton, N.J., 1980.
\bibitem[MW]{milnewaterhouse} J. Milne, W. Waterhouse, Abelian varieties over finite fields, Proc. Sympos. Pure Math., Vol. {\bf XX}, State Univ. New York, Stony Brook, N.Y., 1969, 53--64.
\bibitem[Mu]{mumford} D. Mumford, Abelian varieties, Tata Institute of Fundamental Research Studies in Mathematics, No. {\bf 5},  Published for the Tata Institute of Fundamental Research, Bombay; Oxford University Press, London 1970.
\bibitem[Ne]{neukirch} J. Neukirch, Algebraic number theory, Fundamental Principles of Mathematical Sciences, {\bf 322}, Springer-Verlag, Berlin, 1999.
\bibitem[PR]{pr} V. Platonov, A. Rapinchuk, Algebraic groups and number theory, Pure and Applied Mathematics, {\bf 139}, Academic Press, Inc., Boston, MA, 1994.\bibitem[R1]{ravenel} D. Ravenel, Preprint of part I, available at: http://www.math.rochester.edu/people/faculty/doug/preprints.html
\bibitem[R2]{ravenelbook} D. Ravenel, Nilpotence and periodicity in stable homotopy theory, Annals of Mathematics Studies, {\bf 128}, Princeton University Press, Princeton, NJ, 1992.
\bibitem[Re]{reiner} I. Reiner, Maximal orders, London Mathematical Society Monographs, New Series, {\bf 28}, The Clarendon Press, Oxford University Press, Oxford, 2003.
\bibitem[S]{springer} T. Springer, Linear Algebraic Groups, in: Algebraic geometry  IV, Encyclopaedia of Mathematical Sciences {\bf 55}, Springer-Verlag, Berlin, 1994.
\bibitem[SGA III$_{1}$]{sga3} Sch\'emas en groupes I: Propri\'et\'es g\'en\'erales des sch\'emas en groupes, S\'eminaire de G\'eom\'etrie Alg\'ebrique du Bois Marie 1962/64, Dirig\'e par M. Demazure et A. Grothendieck, Lecture Notes in Mathematics, Vol. {\bf 151}, Springer-Verlag, Berlin-New York 1970.
\bibitem[T]{tatepdiv} J. Tate, $p$-divisible groups, 1967 Proc. Conf. Local Fields (Driebergen, 1966) pp. 158--183.
\bibitem[Wa1]{waterhouse} W. Waterhouse, Abelian varieties over finite fields, Ann. Sci. \'Ecole Norm. Sup. (4) {\bf 2} (1969), 521--560.
\bibitem[Wa2]{waterhousegroups} W. Waterhouse, Introduction to affine group schemes, Graduate Texts in Mathematics, {\bf 66}, Springer-Verlag, New York-Berlin, 1979.
\bibitem[We]{weil} A. Weil, Ad\`eles and algebraic groups, Progress in Mathematics {\bf 23}, Birkh\"auser, Boston, Mass., 1982.

\end{thebibliography}
\end{document}